\documentclass[12pt]{article}
\usepackage{amsfonts}
\usepackage{graphicx}
\usepackage{amsmath}

\begin{document}

\begin{center}
{\Large Lower Bound Estimates of the Order of Meromorphic Solutions to
Non-Homogeneous Linear Differential-Difference Equations }

\quad

\textbf{Rachid BELLAAMA and} \textbf{Benharrat BELA\"{I}DI\footnote{%
Corresponding author} }

\quad

\textbf{Department of Mathematics }

\textbf{Laboratory of Pure and Applied Mathematics }

\textbf{University of Mostaganem (UMAB) }

\textbf{B. P. 227 Mostaganem-(Algeria)}

\textbf{rachidbellaama10@gmail.com}

\textbf{benharrat.belaidi@univ-mosta.dz}

\quad
\end{center}

\noindent \textbf{Abstract. }\ In this article, we deal with the order of
growth of solutions of non-homogeneous linear differential-difference
equation
\begin{equation*}
\sum_{i=0}^{n}\sum_{j=0}^{m}A_{ij}f^{(j)}(z+c_{i})=F(z),
\end{equation*}%
where $A_{ij},$ $F\left( z\right) $ are entire or meromorphic functions and $%
c_{i}$ $\left( 0,1,...,n\right) $ are non-zero distinct complex numbers.
Under the sufficient condition that there exists one coefficient having the
maximal lower order or having the maximal lower type strictly greater than
the order or the type of other coefficients, we obtain estimates of the
lower bound of the order of meromorphic solutions of the above equation.

\quad



\section{Introduction and statement of main results}

\noindent Throughout this article, a meromorphic function means a function
that is meromorphic in the whole complex plane $\mathbb{C}$. We use the
basic notations such as $m\left( r,f\right) ,$ $N\left( r,f\right) ,$ $%
T\left( r,f\right) $ and fundamental results of Nevanlinna's value
distribution theory $(\left[ 8,12,24\right] $). Further, we denote
respectively by $\rho (f),\mu (f),\tau (f),\underline{\tau }(f),$ the order,
the lower order, the type and the lower type of a meromorphic function $f$,
also when $f$ is an entire function, we use the notations $\tau _{M}\left(
f\right) ,\underline{\tau }_{M}(f)$ respectively for the type and lower type
of $f$ (see e.g. $\left[ 8,12,14,24\right] $).

\quad

\noindent Recently, a lot of results were rapidly obtained for complex
difference and complex difference equations $%
([1-3,5,7,10,11,15-17,19,25,26]) $. The back-ground for these studies lies
in the recent difference counterparts of Nevanlinna theory. The key result
here is the difference analogue of the lemma on the logarithmic derivative
obtained by Halburd-Korhonen $\left[ 10,11\right] $ and Chiang-Feng $\left[ 7%
\right] $, independently. Several authors have investigated the properties
of meromorphic solutions of complex linear difference equation%
\begin{equation}
A_{n}f(z+c_{n})+A_{n-1}f(z+c_{n-1})+\cdots
+A_{1}f(z+c_{1})+A_{0}f(z)=A_{n+1},  \tag{1.1}
\end{equation}%
where $A_{j}(z)$ $(j=0,...,n+1)$ are entire or meromorphic functions and $%
c_{i}$ $\left( 1,...,n\right) $ are non-zero distinct complex numbers, when
one coefficient has maximal order (lower order) or among coefficients having
the maximal order (lower order), exactly one has its type (lower type)
strictly greater than others and achieved some important results (see e.g. $%
\left[ 2,7,15,19,25\right] $). Very recently in $\left[ 2\right] $, the
authors have studied the growth of meromorphic solutions of $\left(
1.1\right) $ when one coefficient having maximal lower order or having the
maximal lower type stricly greater than the order or the type of other
coefficients, we have obtained the following theorems.

\quad

\noindent \textbf{Theorem A} $\left( \left[ 2\right] \right) $ \textit{Let }$%
A_{j}(z)$ $(j=0,...,n+1)$\textit{\ be entire functions, and let }$k,l\in
\{0,1,...,n+1\}$\textit{. If the following three assumptions hold
simultaneously}:

\begin{description}
\item[$\left( 1\right) $] $\max \{\mu (A_{k}),\rho (A_{j}),j\neq k,l\}=\rho
\leq \mu (A_{l})<\infty ,\mu (A_{l})>0;$

\item[$\left( 2\right) $] $\underline{\tau }_{M}(A_{l})>\underline{\tau }%
_{M}(A_{k})$, \textit{when} $\mu (A_{l})=\mu (A_{k});$

\item[$\left( 3\right) $] $\max \{\tau _{M}(A_{j}):\rho (A_{j})=\mu
(A_{l}),j\neq k,l\}=\tau _{1}<\underline{\tau }_{M}(A_{l})$, \textit{when} $%
\mu (A_{l})=\max \{\rho (A_{j}),j\neq k,l\}.$
\end{description}

\noindent \textit{Then every meromorphic solution }$f$\textit{\ of }$(1.1)$%
\textit{\ satisfies }$\rho (f)\geq \mu (A_{l})$\textit{\ if} $%
A_{n+1}\not\equiv 0.$ \textit{Furthermore, if }$A_{n+1}(z)\equiv 0,$\textit{%
\ then every meromorphic solution }$f\not\equiv 0$\textit{\ of }$(1.1)$%
\textit{\ satisfies }$\rho (f)\geq \mu (A_{l})+1.$\textit{\ }

\quad

\noindent \textbf{Theorem B} $\left( \left[ 2\right] \right) $ \textit{Let }$%
A_{j}(z)$ $(j=0,...,n+1)$\textit{\ be meromorphic functions, and let }$%
k,l\in \{0,1,...,n+1\}$\textit{. If the following five assumptions hold
simultaneously}:

\begin{description}
\item[$\left( 1\right) $] $\max \{\mu (A_{k}),\rho (A_{j}),j\neq k,l\}=\rho
\leq \mu (A_{l})<\infty ;$

\item[$\left( 2\right) $] $\underline{\tau }(A_{l})>\underline{\tau }(A_{k})$%
, \textit{when} $\mu (A_{l})=\mu (A_{k});$

\item[$\left( 3\right) $]
\begin{equation*}
\tau _{1}=\sum_{\rho (A_{j})=\mu (A_{l}),\,j\neq l,k}\tau (A_{j})<\underline{%
\tau }(A_{l})<+\infty
\end{equation*}%
\textit{when} $\mu (A_{l})=\max \{\rho (A_{j}),j\neq l,k\};$

\item[$\left( 4\right) $]
\begin{equation*}
\tau_{1}+\underline{\tau }(A_{k})<\underline{\tau }(A_{l})<+\infty
\end{equation*}%
\textit{when} $\mu (A_{l})=\mu \left( A_{k}\right) =\max \{\rho(A_{j}),j\neq
k,l\};$

\item[$\left( 5\right) $] $\lambda \left( \frac{1}{A_{l}}\right) <\mu
(A_{l})<\infty .$
\end{description}

\noindent \textit{Then every meromorphic solution }$f$\textit{\ of }$(1.1)$%
\textit{\ satisfies }$\rho (f)\geq \mu (A_{l})$\textit{\ if} $%
A_{n+1}\not\equiv 0.$ \textit{Furthermore, if }$A_{n+1}(z)\equiv 0,$\textit{%
\ then every meromorphic solution }$f\not\equiv 0$\textit{\ of }$(1.1)$%
\textit{\ satisfies }$\rho (f)\geq \mu (A_{l})+1.$\textit{\ }

\quad

\noindent \qquad Historically, the study of complex differential-difference
equations can be traced back to Naftalevich's research. By making use of
operator theory and iteration method, Naftalevich $\left[ 20\right] $
considered the meromorphic solutions on complex differential-difference
equations. However, there are few investigations on complex
differential-difference field using Nevanlinna theory. Recently, some papers
(see, $\left[ 6,18,21,22,26\right] $) focusing on complex
differential-difference emerged. In $\left[ 6\right] ,$ Chen and Zheng have
investigated the growth of solutions of the homogeneous
differential-difference equation
\begin{equation}
\sum_{i=0}^{n}\sum_{j=0}^{m}A_{ij}f^{(j)}(z+c_{i})=0  \tag{1.2}
\end{equation}%
and obtained the following results.

\quad

\noindent \textbf{Theorem C }$\left( \left[ 6\right] \right) $ \textit{Let }$%
A_{ij}\left( z\right) \left( i=0,...,n;j=0,...,m\right) $\textit{\ }\textit{%
\ be entire functions, and}\textit{\ }$a,l\in \{0,1,...,n\},b\in
\{0,1,...,m\}$\textit{\ such that}\textit{\ }$\left( a,b\right) \neq \left(
l,0\right) $\textit{. If the following three assumptions hold simultaneously
}: \textit{\newline
}$\left( \text{1}\right) $\textit{\ }$\max \{\mu (A_{ab}),\rho
(A_{ij}),(i,j)\neq (a,b),(l,0)\}=\rho \leq \mu (A_{l0})<\infty ,\mu
(A_{l0})>0;$ \textit{\newline
}$\left( \text{2}\right) $\textit{\ }$\underline{\tau }_{M}(A_{l0})>%
\underline{\tau }_{M}(A_{ab})$, \textit{when} $\mu {(A_{l0})}=\mu {(A_{ab})}%
; $ \textit{\newline
}$\left( \text{3}\right) $\textit{\ }$\underline{\tau }_{M}(A_{l0})>\max
\{\tau _{M}{(A_{ij})}:\rho (A_{ij})=\mu (A_{l0}),(i,j)\neq (a,b),(l,0)\}$,
when $\mu (A_{l0})=\max \{\rho (A_{ij}):(i,j)\neq (a,b),(l,0)\}.$ \textit{\
Then any non zero meromorphic solution }$f$\textit{\ of }$\left( 1.2\right) $%
\textit{\ satisfies }$\rho (f)\geq \mu (A_{l0})+1$.

\quad

\noindent \textbf{Theorem D} $\left( \left[ 6\right] \right) $ \textit{Let }$%
A_{ij}\left( z\right) $\textit{\ }$\left( i=0,...,n;j=0,...,m\right) $%
\textit{\ be meromorphic functions, and}\textit{\ }$a,l\in
\{0,1,...,n\},b\in \{0,1,...,m\}$\textit{\ such that}\textit{\ }$\left(
a,b\right) \neq \left( l,0\right) .$\textit{\ If the following four
assumptions hold simultaneously}:\textit{\newline
}$\left( \text{1}\right) $\textit{\ }$\delta (\infty ,A_{l0})=\underset{%
r\rightarrow +\infty }{\lim \inf }\frac{m(r,A_{l0})}{T\left( r,A_{l0}\right)
}=\delta >0;$\textit{\newline
}$\left( \text{2}\right) $\textit{\ }$\max \{\mu (A_{ab}),\rho
(A_{ij}),(i,j)\neq (a,b),(l,0)\}=\rho \leq \mu (A_{l0})<\infty ,\mu
(A_{l0})>0;$ \textit{\newline
}$\left( \text{3}\right) $\textit{\ }$\delta \underline{\tau }(A_{l0})>%
\underline{\tau }(A_{ab}),$\textit{\ when}\textit{\ }$\mu {(A_{l0})}=\mu {%
(A_{ab})};$ \textit{\newline
}$\left( \text{4}\right) $\textit{\ }$\delta \underline{\tau }(A_{l0})>\max
\{\tau {(A_{ij})}:\rho (A_{ij})=\mu (A_{l0}),(i,j)\neq (a,b),(l,0)\}$\textit{%
\ when } $\mu (A_{l0})=\max \{\rho {(A_{ij})}:(i,j)\neq (a,b),(l,0)\}.$
\textit{\ Then any non zero meromorphic solution}\textit{\ }$f$\textit{\ of}%
\textit{\ }$\left( 1.2\right) $\textit{\ satisfies}\textit{\ }$\rho (f)\geq
\mu (A_{l0})+1$.

\quad

\noindent \qquad In this paper, by combining complex differential and
difference equations, we extend the results of Theorems C and D for the
complex non-homogeneous linear differential-difference equation
\begin{equation}
\sum_{i=0}^{n}\sum_{j=0}^{m}A_{ij}f^{(j)}(z+c_{i})=F.  \tag{1.3}
\end{equation}%
The main results of this paper state as follows.

\quad

\noindent \textbf{Theorem 1.1} \textit{Let }$%
A_{ij}(z)(i=0,...,n;j=0,...,m),F(z)$\textit{\ be entire functions, and let }$%
k,p$ $(0\leq k\leq n,0\leq p\leq m)$\textit{\ be integers. If there exists
an integer }$l$ $(0\leq l\leq n)$\textit{, such that }$A_{l0}$ \textit{\
satisfies the following three assumptions simultaneously}:

\begin{description}
\item[$\left( 1\right) $] $\max \{\mu (A_{kp}),\rho (F),\rho
(A_{ij}),(i,j)\neq (l,0),(k,p)\}=\rho \leq \mu (A_{l0})<\infty ,$ $\mu
(A_{l0})>0;$

\item[$\left( 2\right) $] $\underline{\tau }_{M}(A_{l0})>\underline{\tau }%
_{M}(A_{kp})$, \textit{when} $\mu (A_{l0})=\mu (A_{kp});$

\item[$\left( 3\right) $] $\max \{\tau _{M}(S):\rho (S)=\mu (A_{l0}),$ $S\in
\{F,A_{ij}:(i,j)\neq (l,0),(k,p)\}\}=\tau _{1}<\underline{\tau }_{M}(A_{l0})$%
, \textit{when} $\mu (A_{l0})=\max \{\rho (F),\rho (A_{ij}),$ $(i,j)\neq
(l,0),(k,p)\}.$
\end{description}

\noindent \textit{Then every meromorphic solution }$f$\textit{\ of }$\left(
1.3\right) $\textit{\ satisfies }$\rho (f)\geq \mu (A_{l0})$\textit{\ if }$%
F(z)\not\equiv 0.$\textit{\ Further, if }$F(z)\equiv 0,$\textit{\ then every
meromorphic solution }$f\not\equiv 0$\textit{\ of }$(1.2)$\textit{\
satisfies }$\rho (f)\geq \mu (A_{l0})+1.$

\quad

\noindent \textbf{Theorem 1.2} \textit{Let }$%
A_{ij}(z)(i=0,...,n;j=0,...,m),F(z)$\textit{\ be meromorphic functions, and
let }$k,p$ $(0\leq k\leq n,0\leq p\leq m)$\ \textit{be integers. If there
exists an integer }$l(0\leq l\leq n)$\textit{, such that }$A_{l0}$ \textit{%
satisfies the following five assumptions simultaneously}:

\begin{description}
\item[$\left( 1\right) $] $\max \{\mu (A_{kp}),\rho (F),\rho
(A_{ij}),(i,j)\neq (l,0),(k,p)\}=\rho \leq \mu (A_{l0})<\infty ;$

\item[$\left( 2\right) $] $\underline{\tau }(A_{l0})>\underline{\tau }%
(A_{kp})$, \textit{when} $\mu (A_{l0})=\mu (A_{kp});$

\item[$\left( 3\right) $]
\begin{equation*}
\tau _{1}=\sum_{\rho (A_{ij})=\mu (A_{l0}),\,(i,j)\neq (l,0),(k,p)}\tau
(A_{ij})+\tau (F)<\underline{\tau }(A_{l0})<+\infty
\end{equation*}%
\textit{when} $\mu (A_{l0})=\max \{\rho (F),\rho (A_{ij}),(i,j)\neq
(l,0),(k,p)\};$

\item[$\left( 4\right) $]
\begin{equation*}
\tau _{1}+\underline{\tau }(A_{kp})<\underline{\tau }(A_{l0})<+\infty
\end{equation*}%
\textit{when} $\mu (A_{l0})=\mu \left( A_{kp}\right) =\max \{\rho (F),\rho
(A_{ij}),(i,j)\neq (l,0),(k,p)\};$

\item[$\left( 5\right) $] $\lambda \left( \frac{1}{A_{l0}}\right) <\mu
(A_{l0})<\infty .$
\end{description}

\noindent \textit{Then every meromorphic solution }$f$\textit{\ of }$\left(
1.3\right) $\textit{\ satisfies }$\rho (f)\geq \mu (A_{l0})$\textit{\ if }$%
F(z)\not\equiv 0.$\textit{\ Further, if }$F(z)\equiv 0,$\textit{\ then every
meromorphic solution }$f\not\equiv 0$\textit{\ of }$(1.2)$\textit{\
satisfies }$\rho (f)\geq \mu (A_{l0})+1.$

\section{Some Preliminary Lemmas}

\noindent \textbf{Lemma 2.1 }$\left( \left[ 9\right] \right) $ \textit{Let }$%
f$\textit{\ be a transcendental meromorphic function of finite order }$\rho
(f)$\textit{, and let k and j be integers satisfying }$k>j\geq 0$\textit{.
Then for any given }$\varepsilon (>0)$\textit{, there exists a subset }$%
E_{1}\subset (1,+\infty )$\textit{\ which has finite logarithmic measure
such that for all }$z$\textit{\ satisfying }$\left\vert z\right\vert
=r\notin \lbrack 0,1]\cup E_{1}$\textit{, we have }%
\begin{equation*}
\left\vert \frac{f^{(k)}(z)}{f^{(j)}(z)}\right\vert \leq |z|^{(k-j)(\rho
(f)-1+\varepsilon )}.
\end{equation*}%
\textbf{Lemma 2.2} $\left( \left[ 7\right] \right) $ \textbf{\ }\textit{Let }%
$f$\textit{\ be a meromorphic function of finite order }$\rho $\textit{, and
let }$c_{1},c_{2}(c_{1}\neq c_{2})$\textit{\ be two arbitrary complex
numbers. Let }$\varepsilon >0$\textit{\ be given, then there exists a subset
}$E_{2}\subset \left( 1,+\infty \right) $\textit{\ with finite logarithmic
measure such that for all }$z$ \textit{satisfying }$\left\vert z\right\vert
=r\notin \lbrack 0,1]\cup E_{2}$\textit{, we have }%
\begin{equation*}
\exp \{-r^{\rho -1+\varepsilon }\}\leq \left\vert \frac{f(z+c_{1})}{%
f(z+c_{2})}\right\vert \leq \exp \{r^{\rho -1+\varepsilon }\}.
\end{equation*}%
\textbf{Lemma 2.3 }$\left( \left[ 8\right] \right) $ \textit{Let }$f$\textit{%
\ be a meromorphic function, }$c$\textit{\ be a non-zero complex constant.
Then we have that for }$r\longrightarrow +\infty $%
\begin{equation*}
(1+o(1))T(r-|c|,f(z))\leq T(r,f(z+c))\leq (1+o(1))T(r+|c|,f(z)).
\end{equation*}%
\textit{Consequently }%
\begin{equation*}
\rho (f(z+c))=\rho (f),\quad \mu (f(z+c))=\mu (f).
\end{equation*}%
\textbf{Lemma 2.4 }$\left( \left[ 4\right] \right) $\textbf{\ }\textit{Let }$%
f$\textit{\ be a meromorphic function of finite order }$\rho $\textit{. Then
for any given }$\varepsilon >0$\textit{, there exists a set }$E_{3}\subset
\left( 1,+\infty \right) $\textit{\ having finite linear measure and finite
logarithmic measure such that for all }$z$\textit{\ satisfying }$\left\vert
z\right\vert =r\notin \lbrack 0,1]\cup E_{3}$ \textit{and sufficiently large
}$r$\textit{, we have }%
\begin{equation*}
\exp \{-r^{\rho +\varepsilon }\}\leq \left\vert f(z)\right\vert \leq \exp
\{r^{\rho +\varepsilon }\}.
\end{equation*}%
\textbf{Lemma 2.5 }$\left( \left[ 13\right] \right) $\textbf{\ }\textit{Let }%
$f$\textit{\ be an entire function with }$\mu (f)<\infty $\textit{. Then for
any given }$\varepsilon (>0)$\textit{, there exists a subset }$E_{4}\subset
(1,+\infty )$\textit{\ with infinite logarithmic measure such that for all }$%
r\in E_{4}$\textit{, we have }%
\begin{equation*}
\mu (f)=\underset{\underset{r\in E_{4}}{r\rightarrow +\infty }}{\lim }\frac{%
\log \log M(r,f)}{\log r}
\end{equation*}%
\textit{and}%
\begin{equation*}
M(r,f)<\exp \{r^{\mu (f)+\varepsilon }\}.
\end{equation*}%
\textbf{Lemma 2.6 }$\left( \left[ 23\right] \right) $\textbf{\ }\textit{Let }%
$f$\textit{\ be an entire function with }$0<\mu (f)<\infty $\textit{. Then
for any given }$\varepsilon (>0)$\textit{, there exists a subset }$%
E_{5}\subset (1,+\infty )$\textit{\ with infinite logarithmic measure such
that for all }$r\in E_{5}$\textit{, we have }%
\begin{equation*}
\underline{\tau }_{M}(f)=\underset{\underset{r\in E_{5}}{r\rightarrow
+\infty }}{\lim }\frac{\log M(r,f)}{\log r}
\end{equation*}%
\textit{and}%
\begin{equation*}
M(r,f)<\exp \{(\underline{\tau }_{M}(f)+\varepsilon )r^{\mu (f)}\}.
\end{equation*}%
\textbf{Lemma 2.7 }$\left( \left[ 7\right] \right) $ \textit{Let }$f$\textit{%
\ be a meromorphic function of finite order }$\rho (f)<\infty $\textit{, and
let }$c_{1},c_{2}$\textit{\ be two distinct complex numbers. Then for each }$%
\varepsilon >0,$\textit{\ we have }%
\begin{equation*}
m\left( r,\frac{f(z+c_{1})}{f(z+c_{2})}\right) =O(r^{\rho (f)-1+\varepsilon
}).
\end{equation*}%
\textbf{Lemma 2.8 }$\left( \left[ 12\right] \right) $ (\textit{Logarithmic
Derivative Lemma}) \textit{Let }$f$\textit{\ be a meromorphic function and }$%
k\geq 1$\textit{\ be an integer. Then, we have}\textit{\ }
\begin{equation*}
m\left( r,\frac{f^{(k)}}{f}\right) =O\left( \log T(r,f)+\log r\right)
\end{equation*}%
\textit{\ possibly outside a set }$E_{6}\subset \lbrack 0,+\infty )$\textit{%
\ of a finite linear measure.} \textit{\ If }$\rho (f)<\infty ,$\textit{\
then}\textit{\ }
\begin{equation*}
m\left( r,\frac{f^{(k)}}{f}\right) =O\left( \log r\right) .
\end{equation*}%
\textbf{Lemma 2.9 }$\left( \left[ 25\right] \right) $ \textit{Let }$f$%
\textit{\ be a meromorphic function with }$\mu (f)<\infty $\textit{. Then
for any given }$\varepsilon (>0)$\textit{, there exists a subset }$%
E_{7}\subset (1,+\infty )$\textit{\ with infinite logarithmic measure such
that for all }$r\in E_{7}$\textit{, we have }%
\begin{equation*}
T(r,f)<r^{\mu (f)+\varepsilon }.
\end{equation*}%
\textbf{Lemma 2.10 }$\left( \left[ 19\right] \right) $ \textit{Let }$f$%
\textit{\ be a meromorphic function with }$0<\mu (f)<\infty $\textit{. Then
for any given }$\varepsilon (>0)$\textit{, there exists a subset }$%
E_{8}\subset (1,+\infty )$\textit{\ with infinite logarithmic measure such
that for all }$r\in E_{8}$\textit{, we have }%
\begin{equation*}
T(r,f)<(\underline{\tau }(f)+\varepsilon )r^{\mu (f)}.
\end{equation*}

\section{Proof of main results}

\noindent \textbf{Proof of Theorem 1.1} If $f$ has infinite order, then the
result holds. Now, we suppose that $\rho (f)<\infty $. \ We divide $(1.3)$
by $f(z+c_{l})$ to get%
\begin{equation*}
-A_{l0}(z)=\sum_{i=0,i\neq l,k}^{n}\sum_{j=0}^{m}A_{ij}\frac{f^{(j)}(z+c_{i})%
}{f(z+c_{i})}\frac{f(z+c_{i})}{f(z+c_{l})}+\sum_{j=0,j\neq p}^{m}A_{kj}\frac{%
f^{(j)}(z+c_{k})}{f(z+c_{k})}\frac{f(z+c_{k})}{f(z+c_{l})}
\end{equation*}%
\begin{equation}
+A_{kp}\frac{f^{(p)}(z+c_{k})}{f(z+c_{k})}\frac{f(z+c_{k})}{f(z+c_{l})}%
+\sum_{j=1}^{m}A_{lj}\frac{f^{(j)}(z+c_{l})}{f(z+c_{l})}-\frac{F(z)}{%
f(z+c_{l})}.  \tag{3.1}
\end{equation}%
Therefore%
\begin{equation*}
|A_{l0}(z)|\leq \sum_{i=0,i\neq l,k}^{n}\sum_{j=0}^{m}|A_{ij}|\left\vert
\frac{f^{(j)}(z+c_{i})}{f(z+c_{i})}\right\vert \left\vert \frac{f(z+c_{i})}{%
f(z+c_{l})}\right\vert
\end{equation*}%
\begin{equation*}
+\sum_{j=0,j\neq p}^{m}|A_{kj}|\left\vert \frac{f^{(j)}(z+c_{k})}{f(z+c_{k})}%
\right\vert \left\vert \frac{f(z+c_{k})}{f(z+c_{l})}\right\vert
\end{equation*}%
\begin{equation}
+|A_{kp}|\left\vert \frac{f^{(p)}(z+c_{k})}{f(z+c_{k})}\right\vert
\left\vert \frac{f(z+c_{k})}{f(z+c_{l})}\right\vert
+\sum_{j=1}^{m}|A_{lj}|\left\vert \frac{f^{(j)}(z+c_{l})}{f(z+c_{l})}%
\right\vert +\left\vert \frac{F(z)}{f(z+c_{l})}\right\vert .  \tag{3.2}
\end{equation}%
From Lemma 2.1 and Lemma 2.3, for any given $\varepsilon (>0)$, there exists
a subset $E_{1}\subset (1,+\infty )$ which has finite logarithmic measure
such that for all $z$ satisfying $\left\vert z\right\vert =r\notin \lbrack
0,1]\cup E_{1}$, we have
\begin{equation}
\left\vert \frac{f^{(j)}(z+c_{i})}{f(z+c_{i})}\right\vert \leq |z|^{j(\rho
\left( f+c_{i}\right) -1+\varepsilon )}=|z|^{j(\rho \left( f\right)
-1+\varepsilon )},\text{ }(i,j)\neq (l,0).  \tag{3.3}
\end{equation}%
It follows by Lemma 2.2 that for any $\varepsilon (>0)$, there exists a
subset $E_{2}\subset \left( 1,+\infty \right) {}$with finite logarithmic
measure such that for all $z$ satisfying $\left\vert z\right\vert =r\notin
\lbrack 0,1]\cup E_{2}$, we have
\begin{equation}
\left\vert \frac{f(z+c_{i})}{f(z+c_{l})}\right\vert \leq \exp \{r^{\rho
\left( f\right) -1+\varepsilon }\},\quad i\neq l.  \tag{3.4}
\end{equation}%
From Lemma 2.3, we get
\begin{equation*}
\rho (f(z+c_{l}))=\rho \left( \frac{1}{f(z+c_{l})}\right) =\rho (f).
\end{equation*}%
So, by Lemma 2.4, for any given $\varepsilon >0$, there exists a subset $%
E_{3}\subset \left( 1,+\infty \right) $ having finite linear measure and
finite logarithmic measure such that for all $z$ satisfying $|z|=r\notin
\lbrack 0,1]\cup E_{3}$ sufficiently large, we have
\begin{equation}
\left\vert \frac{1}{f(z+c_{l})}\right\vert \leq \exp \{r^{\rho
(f)+\varepsilon }\}.  \tag{3.5}
\end{equation}%
In the following, we divide the proof into four cases:

\quad

\noindent \textbf{Case }$\left( \text{i}\right) .$ We suppose that $\rho
<\mu (A_{l0}).$

\quad

\noindent For $S\in \left\{ F,A_{ij}:(i,j)\neq (l,0),(k,p)\right\} ,$ by the
definition of $\rho (S),$ for any given $\varepsilon >0$ and sufficiently
large $r$, we have
\begin{equation}
\left\vert S\left( z\right) \right\vert \leq \exp \{r^{\rho (S)+\varepsilon
}\}\leq \exp \{r^{\rho +\varepsilon }\}.  \tag{3.6}
\end{equation}%
It follows by the definition of $\mu (A_{l0})$, for sufficiently small $%
\varepsilon >0$ and sufficiently large $r$, we have%
\begin{equation}
\left\vert A_{l0}\left( z\right) \right\vert \geq \exp \{r^{\mu
(A_{l0})-\varepsilon }\}.  \tag{3.7}
\end{equation}%
It also follows by the definition of $\mu (A_{kp})$ and Lemma 2.5, for any
given $\varepsilon (>0)$, there exists a subset $E_{4}\subset (1,+\infty )$
with infinite logarithmic measure such that for all $r\in E_{4}$, we have
\begin{equation}
\left\vert A_{kp}\left( z\right) \right\vert \leq \exp \{r^{\mu
(A_{kp})+\varepsilon }\}.  \tag{3.8}
\end{equation}%
By substituting $(3.3)-(3.8)$ into $\left( 3.2\right) ,$ for all $z$
satisfying $|z|=r\in E_{4}\setminus \left( \lbrack 0,1]\cup E_{1}\cup
E_{2}\cup E_{3}\right) $, we obtain%
\begin{equation*}
\exp \{r^{\mu (A_{l0})-\varepsilon }\}\leq \sum_{i=0,i\neq
l,k}^{n}\sum_{j=0}^{m}\exp \{r^{\rho +\varepsilon }\}|z|^{j(\rho
(f)-1+\varepsilon )}\exp \{r^{\rho (f)-1+\varepsilon }\}
\end{equation*}%
\begin{equation*}
+\sum_{j=0,j\neq p}^{m}\exp \{r^{\rho +\varepsilon }\}|z|^{j(\rho
(f)-1+\varepsilon )}\exp \{r^{\rho (f)-1+\varepsilon }\}
\end{equation*}%
\begin{equation*}
+|z|^{p(\rho (f)-1+\varepsilon )}\exp \{r^{\mu (A_{kp})+\varepsilon }\}\exp
\{r^{\rho (f)-1+\varepsilon }\}
\end{equation*}%
\begin{equation*}
+\sum_{j=1}^{m}\exp \{r^{\rho +\varepsilon }\}|z|^{j(\rho (f)-1+\varepsilon
)}+\exp \{r^{\rho +\varepsilon }\}\exp \{r^{\rho (f)+\varepsilon }\}.
\end{equation*}%
\begin{equation*}
\leq \left( (n-1)(m+1)+2m\right) r^{m(\rho (f)-1+\varepsilon )}\exp
\{r^{\rho +\varepsilon }\}\exp \{r^{\rho (f)-1+\varepsilon }\}
\end{equation*}%
\begin{equation}
+r^{p(\rho (f)-1+\varepsilon )}\exp \{r^{\mu (A_{kp})+\varepsilon }\}\exp
\{r^{\rho (f)-1+\varepsilon }\}+\exp \{r^{\rho +\varepsilon }\}\exp
\{r^{\rho (f)+\varepsilon }\}.  \tag{3.9}
\end{equation}%
Now, we may choose sufficiently small $\varepsilon $ satisfying $%
0<3\varepsilon <\mu (A_{l0})-\rho ,$ we deduce from $\left( 3.9\right) $
that for $|z|=r\in E_{4}\setminus ([0,1]\cup E_{1}\cup E_{2}\cup E_{3}),$ $%
r\rightarrow +\infty $
\begin{equation*}
\exp \{r^{\mu (A_{l0})-2\varepsilon }\}\leq \exp \{r^{\rho (f)+\varepsilon
}\}.
\end{equation*}%
Therefore, $\mu (A_{l0})\leq \rho (f)+3\varepsilon $, since $\varepsilon >0$
is arbitrary, then $\rho (f)\geq \mu (A_{l0}).$

\quad

\noindent Further, if \textit{\ }$F\equiv 0,$\textit{\ }then by substituting
$(3.3),$ $(3.4)$ and $(3.6)-(3.8)$ into $\left( 3.2\right) ,$ for all $z$
satisfying $|z|=r\in E_{4}\setminus \left( \lbrack 0,1]\cup E_{1}\cup
E_{2}\right) $, we obtain%
\begin{equation*}
\exp \{r^{\mu (A_{l0})-\varepsilon }\}\leq \left( nm+n+m-1\right) r^{m(\rho
(f)-1+\varepsilon )}\exp \{r^{\rho +\varepsilon }\}\exp \{r^{\rho
(f)-1+\varepsilon }\}
\end{equation*}%
\begin{equation}
+r^{p(\rho (f)-1+\varepsilon )}\exp \{r^{\mu (A_{kp})+\varepsilon }\}\exp
\{r^{\rho (f)-1+\varepsilon }\}.  \tag{3.10}
\end{equation}%
By choosing sufficiently small $\varepsilon $ satisfying $0<3\varepsilon
<\mu (A_{l0})-\rho ,$ we deduce from $\left( 3.10\right) $ that for $%
|z|=r\in E_{4}\setminus ([0,1]\cup E_{1}\cup E_{2}),$ $r\rightarrow +\infty $
\begin{equation*}
\exp \{r^{\mu (A_{l0})-2\varepsilon }\}\leq \exp \{r^{\rho (f)-1+\varepsilon
}\},
\end{equation*}%
that is, $\mu (A_{l0})\leq \rho (f)-1+3\varepsilon $, since $\varepsilon >0$
is arbitrary, then $\rho (f)\geq \mu (A_{l0})+1.$

\noindent \textbf{Case }$\left( \text{ii}\right) $ We suppose that $\beta
=\max \{\rho (F),\rho (A_{ij}),(i,j)\neq (l,0),(k,p)\}<\mu (A_{l0})=\mu
(A_{kp}),$ $\underline{\tau }_{M}(A_{l0})>\underline{\tau }_{M}(A_{kp}).$

\quad

\noindent For $S\in \left\{ F,A_{ij}:(i,j)\neq (l,0),(k,p)\right\} ,$ by the
definition of $\rho (S),$ for any given $\varepsilon (>0),$ and sufficiently
large $r$, we have
\begin{equation}
\left\vert S\left( z\right) \right\vert \leq \exp \{r^{\rho (S)+\varepsilon
}\}\leq \exp \{r^{\beta +\varepsilon }\}.  \tag{3.11}
\end{equation}%
From the definition of $\underline{\tau }_{M}(A_{l0})$, for sufficiently
small $\varepsilon >0$ and sufficiently large $r$, we have
\begin{equation}
\left\vert A_{l0}\left( z\right) \right\vert \geq \exp \{\left( \underline{%
\tau }_{M}(A_{l0})-\varepsilon \right) r^{\mu (A_{l0})}\}.  \tag{3.12}
\end{equation}%
Also, from the definition of $\underline{\tau }_{M}(A_{kp})$ and Lemma 2.6,
for any given $\varepsilon (>0)$, there exists a subset $E_{5}\subset
(1,+\infty )$ with infinite logarithmic measure such that for all $r\in
E_{5} $, we have
\begin{equation}
\left\vert A_{kp}\left( z\right) \right\vert \leq \exp \{\left( \underline{%
\tau }_{M}(A_{kp})+\varepsilon \right) r^{\mu (A_{kp})}\}=\exp \{\left(
\underline{\tau }_{M}(A_{kp})+\varepsilon \right) r^{\mu (A_{l0})}\}.
\tag{3.13}
\end{equation}%
By substituting $(3.3)-(3.5)$ and $(3.11)-(3.13)$ into $\left( 3.2\right) $,
for all $z$ satisfying $|z|=r\in E_{5}\setminus ([0,1]\cup E_{1}\cup
E_{2}\cup E_{3})$, we get%
\begin{equation*}
\exp \{(\underline{\tau }_{M}(A_{l0})-\varepsilon )r^{\mu (A_{l0})}\}
\end{equation*}%
\begin{equation*}
\leq \left( nm+n+m-1\right) r^{m(\rho (f)-1+\varepsilon )}\exp \{r^{\beta
+\varepsilon }\}\exp \{r^{\rho (f)-1+\varepsilon }\}
\end{equation*}%
\begin{equation}
+r^{p(\rho (f)-1+\varepsilon )}\exp \{(\underline{\tau }_{M}(A_{kp})+%
\varepsilon )r^{\mu (A_{l0})}\}\exp \{r^{\rho (f)-1+\varepsilon }\}+\exp
\{r^{\beta +\varepsilon }\}\exp \{r^{\rho (f)+\varepsilon }\}.  \tag{3.14}
\end{equation}%
Therefore, we may choose sufficiently small $\varepsilon $, $0<2\varepsilon
<\min \{\mu (A_{l0})-\beta ,\underline{\tau }_{M}(A_{l0})-\underline{\tau }%
_{M}(A_{kp})\}$, then from $(3.14)$ for $r\in E_{5}\setminus ([0,1]\cup
E_{1}\cup E_{2}\cup E_{3})$ sufficiently large, we obtain
\begin{equation*}
\exp \{(\underline{\tau }_{M}(A_{l0})-\underline{\tau }_{M}(A_{kp})-2%
\varepsilon )r^{\mu (A_{l0})-\varepsilon }\}\leq \exp \{r^{\rho
(f)+\varepsilon }\}.
\end{equation*}%
Then, $\mu (A_{l0})\leq \rho (f)+2\varepsilon $, since $\varepsilon >0$ is
arbitrary, so $\rho (f)\geq \mu (A_{l0}).$

\quad

\noindent Further, if \textit{\ }$F\equiv 0,$\textit{\ }then by substituting
$(3.3),$ $(3.4)$ and $(3.11)-(3.13)$ into $\left( 3.2\right) $, for all $z$
satisfying $|z|=r\in E_{5}\setminus ([0,1]\cup E_{1}\cup E_{2})$, we have%
\begin{equation*}
\exp \{(\underline{\tau }_{M}(A_{l0})-\varepsilon )r^{\mu (A_{l0})}\}
\end{equation*}%
\begin{equation*}
\leq \left( nm+n+m-1\right) r^{m(\rho (f)-1+\varepsilon )}\exp \{r^{\beta
+\varepsilon }\}\exp \{r^{\rho (f)-1+\varepsilon }\}
\end{equation*}%
\begin{equation}
+r^{p(\rho (f)-1+\varepsilon )}\exp \{(\underline{\tau }_{M}(A_{kp})+%
\varepsilon )r^{\mu (A_{l0})}\}\exp \{r^{\rho (f)-1+\varepsilon }\}.
\tag{3.15}
\end{equation}%
Now, we may choose sufficiently small $\varepsilon $, $0<2\varepsilon <\min
\{\mu (A_{l0})-\beta ,\underline{\tau }_{M}(A_{l0})-\underline{\tau }%
_{M}(A_{kp})\}$, then from $(3.15)$ for $r\in E_{5}\setminus ([0,1]\cup
E_{1}\cup E_{2})$ sufficiently large, we get
\begin{equation*}
\exp \{(\underline{\tau }_{M}(A_{l0})-\underline{\tau }_{M}(A_{kp})-2%
\varepsilon )r^{\mu (A_{l0})-\varepsilon }\}\leq \exp \{r^{\rho
(f)-1+\varepsilon }\},
\end{equation*}%
that is, $\mu (A_{l0})\leq \rho (f)-1+2\varepsilon $, since $\varepsilon >0$
is arbitrary, then $\rho (f)\geq \mu (A_{l0})+1.$

\quad

\noindent \textbf{Case }$\left( \text{iii}\right) $ When $\mu (A_{l0})=\max
\{\rho (F),\rho (A_{ij}),$ $(i,j)\neq (l,0),(k,p)\}>\mu (A_{kp})$, and $\max
\{\tau _{M}(S):\rho (S)=\mu (A_{l0}),S\in \{F,A_{ij}:(i,j)\neq
(l,0),(k,p)\}\}=\tau _{1}<\underline{\tau }_{M}(A_{l0}).$

\quad

\noindent For $S\in \left\{ F,A_{ij}:(i,j)\neq (l,0),(k,p)\right\} ,$ by the
definitions of $\rho (S)$ and $\tau _{M}(S)$, for any given $\varepsilon >0$
and sufficiently large $r$, we have%
\begin{equation}
\left\vert S\left( z\right) \right\vert \leq \left\{
\begin{array}{c}
\exp \{r^{\rho (S)+\varepsilon }\}\leq \exp \{r^{\mu (A_{l0})-\varepsilon
}\},\text{ if }\rho (S)<\mu (A_{l0}), \\
\exp \{\left( \tau _{1}+\varepsilon \right) r^{\mu (A_{l0})}\},\text{ if }%
\rho (S)=\mu (A_{l0}).%
\end{array}%
\right.  \tag{3.16}
\end{equation}%
Then, by substituting $\left( 3.3\right) -\left( 3.5\right) ,$ $\left(
3.8\right) ,$ $\left( 3.12\right) $ and $\left( 3.16\right) $ into $\left(
3.2\right) $, for all $z$ satisfying $\left\vert z\right\vert =r\in
E_{4}\setminus ([0,1]\cup E_{1}\cup E_{2}\cup E_{3})$ sufficiently large, we
obtain%
\begin{equation*}
\exp \{\left( \underline{\tau }_{M}(A_{l0})-\varepsilon \right) r^{\mu
(A_{l0})}\}
\end{equation*}%
\begin{equation*}
\leq O\left( r^{m(\rho (f)-1+\varepsilon )}\exp \{\left( \tau
_{1}+\varepsilon \right) r^{\mu (A_{l0})}\}\exp \left\{ r^{\rho
(f)-1+\varepsilon }\right\} \right)
\end{equation*}%
\begin{equation*}
+O\left( r^{m(\rho (f)-1+\varepsilon )}\exp \{r^{\mu (A_{l0})-\varepsilon
}\}\exp \left\{ r^{\rho (f)-1+\varepsilon }\right\} \right)
\end{equation*}%
\begin{equation*}
+r^{p(\rho (f)-1+\varepsilon )}\exp \{r^{\mu (A_{kp})+\varepsilon }\}\exp
\left\{ r^{\rho (f)-1+\varepsilon }\right\}
\end{equation*}%
\begin{equation*}
+O\left( r^{m(\rho (f)-1+\varepsilon )}\exp \{r^{\mu (A_{l0})-\varepsilon
}\}\right)
\end{equation*}%
\begin{equation}
+O\left( r^{m(\rho (f)-1+\varepsilon )}\exp \left\{ \left( \tau
_{1}+\varepsilon \right) r^{\mu (A_{l0})}\right\} \right) +\exp \{\left(
\tau _{1}+\varepsilon \right) r^{\mu (A_{l0})}\}\exp \left\{ r^{\rho
(f)+\varepsilon }\right\} .  \tag{3.17}
\end{equation}%
Now, we may choose sufficiently small $\varepsilon $ satisfying
\begin{equation*}
0<2\varepsilon <\min \{\mu (A_{l0})-\mu (A_{kp}),\underline{\tau }%
_{M}(A_{l0})-\tau _{1}\},
\end{equation*}%
then from $\left( 3.17\right) $ for sufficiently large $r\in E_{4}\setminus
([0,1]\cup E_{1}\cup E_{2}\cup E_{3})$, we get
\begin{equation*}
\exp \{\left( \underline{\tau }_{M}(A_{l0})-\tau _{1}-2\varepsilon \right)
r^{\mu (A_{l0})-\varepsilon }\}\leq \exp \left\{ r^{\rho (f)+\varepsilon
}\right\} .
\end{equation*}%
That means, $\mu (A_{l0})\leq \rho (f)+2\varepsilon $, since $\varepsilon >0$
is arbitrary, then $\rho (f)\geq \mu (A_{l0}).$

\quad

\noindent Further, if \textit{\ }$F\equiv 0,$\textit{\ }then by substituting
$\left( 3.3\right) ,$ $\left( 3.4\right) ,$ $\left( 3.8\right) ,$ $\left(
3.12\right) $ and $\left( 3.16\right) $ into $\left( 3.2\right) $, for all $%
z $ satisfying $\left\vert z\right\vert =r\in E_{4}\setminus ([0,1]\cup
E_{1}\cup E_{2})$ sufficiently large, we have%
\begin{equation*}
\exp \{\left( \underline{\tau }_{M}(A_{l0})-\varepsilon \right) r^{\mu
(A_{l0})}\}
\end{equation*}%
\begin{equation*}
\leq O\left( r^{m(\rho (f)-1+\varepsilon )}\exp \{\left( \tau
_{1}+\varepsilon \right) r^{\mu (A_{l0})}\}\exp \left\{ r^{\rho
(f)-1+\varepsilon }\right\} \right)
\end{equation*}%
\begin{equation*}
+O\left( r^{m(\rho (f)-1+\varepsilon )}\exp \{r^{\mu (A_{l0})-\varepsilon
}\}\exp \left\{ r^{\rho (f)-1+\varepsilon }\right\} \right)
\end{equation*}%
\begin{equation*}
+r^{p(\rho (f)-1+\varepsilon )}\exp \{r^{\mu (A_{kp})+\varepsilon }\}\exp
\left\{ r^{\rho (f)-1+\varepsilon }\right\}
\end{equation*}%
\begin{equation}
+O\left( r^{m(\rho (f)-1+\varepsilon )}\exp \{r^{\mu (A_{l0})-\varepsilon
}\}\right) +O\left( r^{m(\rho (f)-1+\varepsilon )}\exp \left\{ \left( \tau
_{1}+\varepsilon \right) r^{\mu (A_{l0})}\right\} \right) .  \tag{3.18}
\end{equation}%
Now, we may choose sufficiently small $\varepsilon $ satisfying
\begin{equation*}
0<2\varepsilon <\min \{\mu (A_{l0})-\mu (A_{kp}),\underline{\tau }%
_{M}(A_{l0})-\tau _{1}\},
\end{equation*}%
then from $\left( 3.18\right) $ for sufficiently large $r\in E_{4}\setminus
([0,1]\cup E_{1}\cup E_{2})$, we get
\begin{equation*}
\exp \{\left( \underline{\tau }_{M}(A_{l0})-\tau _{1}-2\varepsilon \right)
r^{\mu (A_{l0})-\varepsilon }\}\leq \exp \left\{ r^{\rho (f)-1+\varepsilon
}\right\} .
\end{equation*}%
That means, $\mu (A_{l0})\leq \rho (f)-1+2\varepsilon $, since $\varepsilon
>0$ is arbitrary, then $\rho (f)\geq \mu (A_{l0})+1.$

\quad

\noindent \textbf{Case} $\left( \text{iv}\right) $ We suppose that $\max
\{\rho (A_{ij}),\rho (F),(i,j)\neq (l,0),(k,p)\}=\mu (A_{kp})=\mu (A_{l0})$
and $\max \{\underline{\tau }_{M}(A_{kp}),\tau _{M}(S):\rho (S)=\mu
(A_{l0}),S\in \{F,A_{ij}:(i,j)\neq (l,0),(k,p)\}\}=\tau _{2}<\underline{\tau
}_{M}(A_{l0}).$

\quad

\noindent It follows by substituting $\left( 3.3\right) -\left( 3.5\right) ,$
$\left( 3.12\right) ,$ $\left( 3.13\right) $ and $\left( 3.16\right) $ into $%
\left( 3.2\right) $, for all $z$ satisfying $\left\vert z\right\vert =r\in
E_{5}\setminus ([0,1]\cup E_{1}\cup E_{2}\cup E_{3})$ sufficiently large, we
have%
\begin{equation*}
\exp \{\left( \underline{\tau }_{M}(A_{l0})-\varepsilon \right) r^{\mu
(A_{l0})}\}
\end{equation*}%
\begin{equation*}
\leq O\left( r^{m(\rho (f)-1+\varepsilon )}\exp \{\left( \tau
_{2}+\varepsilon \right) r^{\mu (A_{l0})}\}\exp \left\{ r^{\rho
(f)-1+\varepsilon }\right\} \right)
\end{equation*}%
\begin{equation*}
+O\left( r^{m(\rho (f)-1+\varepsilon )}\exp \{r^{\mu (A_{l0})-\varepsilon
}\}\exp \left\{ r^{\rho (f)-1+\varepsilon }\right\} \right)
\end{equation*}%
\begin{equation*}
+r^{p(\rho (f)-1+\varepsilon )}\exp \{\left( \underline{\tau }%
_{M}(A_{kp})+\varepsilon \right) r^{\mu (A_{l0})}\}\exp \left\{ r^{\rho
(f)-1+\varepsilon }\right\}
\end{equation*}%
\begin{equation*}
+O\left( r^{m(\rho (f)-1+\varepsilon )}\exp \{r^{\mu (A_{l0})-\varepsilon
}\}\right) +O\left( r^{m(\rho (f)-1+\varepsilon )}\exp \left\{ \left( \tau
_{2}+\varepsilon \right) r^{\mu (A_{l0})}\right\} \right)
\end{equation*}%
\begin{equation}
+\exp \{\left( \tau _{2}+\varepsilon \right) r^{\mu (A_{l0})}\}\exp \left\{
r^{\rho (f)+\varepsilon }\right\} .  \tag{3.19}
\end{equation}%
Now, we may choose sufficiently small $\varepsilon $ satisfying
\begin{equation*}
0<2\varepsilon <\underline{\tau }_{M}(A_{l0})-\tau _{2},
\end{equation*}%
from $\left( 3.19\right) $ for sufficiently large $r\in E_{5}\setminus
([0,1]\cup E_{1}\cup E_{2}\cup E_{3})$, we get
\begin{equation*}
\exp \{\left( \underline{\tau }_{M}(A_{l0})-\tau _{2}-2\varepsilon \right)
r^{\mu (A_{l0})-\varepsilon }\}\leq \exp \left\{ r^{\rho (f)+\varepsilon
}\right\} ,
\end{equation*}%
that means, $\mu (A_{l0})\leq \rho (f)+2\varepsilon ,$ since $\varepsilon >0$
is arbitrary, then $\rho (f)\geq \mu (A_{l0}).$

\quad

\noindent Further, if \textit{\ }$F\equiv 0,$\textit{\ } by substituting $%
\left( 3.3\right) ,$ $\left( 3.4\right) ,$ $\left( 3.12\right) ,$ $\left(
3.13\right) $ and $\left( 3.16\right) $ into $\left( 3.2\right) $, for all $%
z $ satisfying $\left\vert z\right\vert =r\in E_{5}\setminus ([0,1]\cup
E_{1}\cup E_{2})$ sufficiently large, we have%
\begin{equation*}
\exp \{\left( \underline{\tau }_{M}(A_{l0})-\varepsilon \right) r^{\mu
(A_{l0})}\}\leq O\left( r^{m(\rho (f)-1+\varepsilon )}\exp \{\left( \tau
_{2}+\varepsilon \right) r^{\mu (A_{l0})}\}\exp \left\{ r^{\rho
(f)-1+\varepsilon }\right\} \right)
\end{equation*}%
\begin{equation*}
+O\left( r^{m(\rho (f)-1+\varepsilon )}\exp \{r^{\mu (A_{l0})-\varepsilon
}\}\exp \left\{ r^{\rho (f)-1+\varepsilon }\right\} \right)
\end{equation*}%
\begin{equation*}
+r^{p(\rho (f)-1+\varepsilon )}\exp \{\left( \underline{\tau }%
_{M}(A_{kp})+\varepsilon \right) r^{\mu (A_{l0})}\}\exp \left\{ r^{\rho
(f)-1+\varepsilon }\right\}
\end{equation*}%
\begin{equation}
+O\left( r^{m(\rho (f)-1+\varepsilon )}\exp \{r^{\mu (A_{l0})-\varepsilon
}\}\right) +O\left( r^{m(\rho (f)-1+\varepsilon )}\exp \left\{ \left( \tau
_{2}+\varepsilon \right) r^{\mu (A_{l0})}\right\} \right) .  \tag{3.20}
\end{equation}%
Now, we may choose sufficiently small $\varepsilon $ satisfying
\begin{equation*}
0<2\varepsilon <\underline{\tau }_{M}(A_{l0})-\tau _{2},
\end{equation*}%
from $\left( 3.20\right) $ for sufficiently large $r\in E_{5}\setminus
([0,1]\cup E_{1}\cup E_{2})$, we get
\begin{equation*}
\exp \{\left( \underline{\tau }_{M}(A_{l0})-\tau _{2}-2\varepsilon \right)
r^{\mu (A_{l0})-\varepsilon }\}\leq \exp \left\{ r^{\rho (f)-1+\varepsilon
}\right\} .
\end{equation*}%
That means, $\mu (A_{l0})\leq \rho (f)-1+2\varepsilon $, since $\varepsilon
>0$ is arbitrary, then $\rho (f)\geq \mu (A_{l0})+1.$ The proof of Theorem
1.1 is complete.

\quad

\noindent \textbf{Proof of Theorem 1.2 }If $f$ has infinite order, then the
result holds. Now, we suppose that $\rho (f)<\infty $. By $\left( 3.1\right)
,$ we have%
\begin{equation*}
T(r,A_{l0}(z))=m(r,A_{l0}(z))+N(r,A_{l0}(z))
\end{equation*}%
\begin{equation*}
\leq \sum_{i=0,i\neq l,k}^{n}\sum_{j=0}^{m}m(r,A_{ij}(z))+m(r,A_{kp}(z))
\end{equation*}%
\begin{equation*}
+\sum_{j=0,j\neq
p}^{m}m(r,A_{kj}(z))+\sum_{j=1}^{m}m(r,A_{lj}(z))+\sum_{i=0,i\neq
l,k}^{n}\sum_{j=0}^{m}m\left( r,\frac{f^{(j)}(z+c_{i})}{f(z+c_{i})}\right)
\end{equation*}%
\begin{equation*}
+\sum_{i=0,i\neq l,k}^{n}m\left( r,\frac{f(z+c_{i})}{f(z+c_{l})}\right)
+\sum_{j=1}^{m}m\left( r,\frac{f^{(j)}(z+c_{k})}{f(z+c_{k})}\right)
+2m\left( r,\frac{f(z+c_{k})}{f(z+c_{l})}\right)
\end{equation*}%
\begin{equation*}
+\sum_{j=1}^{m}m\left( r,\frac{f^{(j)}(z+c_{l})}{f(z+c_{l})}\right) +m\left(
r,F(z)\right) +m\left( r,\frac{1}{f(z+c_{l})}\right)
\end{equation*}%
\begin{equation*}
+N(r,A_{l0}(z))+O(1)\leq \sum_{i=0,i\neq
l,k}^{n}\sum_{j=0}^{m}T(r,A_{ij}(z))+T(r,A_{kp}(z))
\end{equation*}%
\begin{equation*}
+\sum_{j=0,j\neq
p}^{m}T(r,A_{kj}(z))+\sum_{j=1}^{m}T(r,A_{lj}(z))+\sum_{i=0,i\neq
l,k}^{n}\sum_{j=1}^{m}m\left( r,\frac{f^{(j)}(z+c_{i})}{f(z+c_{i})}\right)
\end{equation*}%
\begin{equation*}
+\sum_{i=0,i\neq l,k}^{n}m\left( r,\frac{f(z+c_{i})}{f(z+c_{l})}\right)
+\sum_{j=1}^{m}m\left( r,\frac{f^{(j)}(z+c_{k})}{f(z+c_{k})}\right)
+2m\left( r,\frac{f(z+c_{k})}{f(z+c_{l})}\right)
\end{equation*}%
\begin{equation*}
+\sum_{j=1}^{m}m\left( r,\frac{f^{(j)}(z+c_{l})}{f(z+c_{l})}\right) +T\left(
r,F(z)\right) +T\left( r,\frac{1}{f(z+c_{l})}\right) +N(r,A_{l0}(z))+O(1).
\end{equation*}%
By Lemma 2.3 and the first main theorem of Nevanlinna, when $r$ sufficiently
large, we have
\begin{equation*}
T\left( r,\frac{1}{f(z+c_{l})}\right) =T\left( r,f(z+c_{l})\right) +O\left(
1\right) \leq (1+o(1))T(r+\left\vert c_{l}\right\vert ,f)\leq 2T(2r,f).
\end{equation*}%
So, for $r$ sufficiently large we obtain%
\begin{equation*}
T(r,A_{l0}(z))\leq \sum_{i=0,i\neq
l,k}^{n}\sum_{j=0}^{m}T(r,A_{ij}(z))+T(r,A_{kp}(z))
\end{equation*}%
\begin{equation*}
+\sum_{j=0,j\neq p}^{m}T(r,A_{kj}(z))+\sum_{j=1}^{m}T(r,A_{lj}(z))
\end{equation*}%
\begin{equation*}
+\sum_{i=0,i\neq l,k}^{n}\sum_{j=1}^{m}m\left( r,\frac{f^{(j)}(z+c_{i})}{%
f(z+c_{i})}\right) +\sum_{i=0,i\neq l,k}^{n}m\left( r,\frac{f(z+c_{i})}{%
f(z+c_{l})}\right) +T\left( r,F(z)\right)
\end{equation*}%
\begin{equation*}
+2T(2r,f)+\sum_{j=1}^{m}m\left( r,\frac{f^{(j)}(z+c_{l})}{f(z+c_{l})}\right)
\end{equation*}%
\begin{equation}
+\sum_{j=1}^{m}m\left( r,\frac{f^{(j)}(z+c_{k})}{f(z+c_{k})}\right)
+2m\left( r,\frac{f(z+c_{k})}{f(z+c_{l})}\right) +N(r,A_{l0}(z))+O(1).
\tag{3.21}
\end{equation}%
By Lemma 2.7, for any given $\varepsilon (>0)$, we have
\begin{equation}
m\left( r,\frac{f(z)}{f(z+c_{l})}\right) =O(r^{\rho (f)-1+\varepsilon
}),\quad m\left( r,\frac{f(z+c_{j})}{f(z+c_{l})}\right) =O(r^{\rho
(f)-1+\varepsilon }),\quad j\neq l.  \tag{3.22}
\end{equation}%
It follows by Lemma 2.8, there exists a subset $E_{6}\subset \lbrack
0,+\infty \lbrack $ of a finite linear measure such that for all $r\notin
E_{6}$ sufficiently large, we have
\begin{equation}
m\left( r,\frac{f^{(j)}(z+c_{i})}{f(z+c_{i})}\right) =O\left( \log r\right)
\quad \left( i=0,...,n;j=1,...,m\right) .  \tag{3.23}
\end{equation}%
From the definition of $\lambda \left( \frac{1}{A_{l0}}\right) $, for any
given $\varepsilon >0$ and sufficiently large $r$, we have
\begin{equation}
N\left( r,A_{l0}\right) \leq r^{\lambda \left( \frac{1}{A_{l0}}\right)
+\varepsilon }.  \tag{3.24}
\end{equation}%
In the following, we divide the proof into four cases:

\quad

\noindent \textbf{Case }$\left( \text{i}\right) $ We suppose that $\rho <\mu
(A_{l0}).$\newline
For $S\in \left\{ F,A_{ij}:(i,j)\neq (l,0),(k,p)\right\} $, from the
definition of $\rho (S)$ and $\rho (f)$ for any given $\varepsilon >0$ and
sufficiently large $r$, we have
\begin{equation}
T(r,S)\leq r^{\rho (S)+\varepsilon }\leq r^{\rho +\varepsilon },  \tag{3.25}
\end{equation}%
\begin{equation}
T(r,f)\leq r^{\rho (f)+\varepsilon }.  \tag{3.26}
\end{equation}%
It follows from the definition of $\mu (A_{l0})$, for sufficiently small $%
\varepsilon >0$ and sufficiently large $r$, we have%
\begin{equation}
T(r,A_{l0})\geq r^{\mu (A_{l0})-\varepsilon }.  \tag{3.27}
\end{equation}%
It follows from the definition of $\mu (A_{kp})$ and Lemma 2.9, for any
given $\varepsilon (>0)$, there exists a subset $E_{7}\subset (1,+\infty )$
with infinite logarithmic measure such that for all $r\in E_{7}$, we have
\begin{equation}
T(r,A_{kp})\leq r^{\mu (A_{kp})+\varepsilon }.  \tag{3.28}
\end{equation}%
By substituting $(3.22)-(3.28)$ into $\left( 3.21\right) $ for sufficiently
large $r\in E_{7}\backslash E_{6}$, we obtain%
\begin{equation*}
r^{\mu (A_{l0})-\varepsilon }\leq \left( (n-1)(m+1)+2m\right) r^{\rho
+\varepsilon }+r^{\mu (A_{kp})+\varepsilon }+O(r^{\rho (f)-1+\varepsilon })
\end{equation*}%
\begin{equation}
+2\left( 2r\right) ^{\rho (f)+\varepsilon }+r^{\rho +\varepsilon
}+r^{\lambda \left( \frac{1}{A_{l0}}\right) +\varepsilon }+O(\log r).
\tag{3.29}
\end{equation}%
We may choose sufficiently small $\varepsilon $ satisfying
\begin{equation*}
0<3\varepsilon <\min \left\{ \mu (A_{l0})-\rho ,\mu (A_{l0})-\lambda \left(
\frac{1}{A_{l0}}\right) \right\} ,
\end{equation*}%
it follows from $\left( 3.29\right) $ that for $r\in E_{7}\backslash E_{6}$,
$r\rightarrow +\infty $
\begin{equation*}
r^{\mu (A_{l0})-2\varepsilon }\leq r^{\rho (f)+\varepsilon },
\end{equation*}%
that means, $\mu (A_{l0})\leq \rho (f)+3\varepsilon $, since $\varepsilon >0$
is arbitrary, then $\rho (f)\geq \mu (A_{l0}).$

\noindent Further, if \textit{\ }$F\equiv 0,$\textit{\ }then by\textit{\ }%
substituting $(3.22)-(3.25),$ $(3.27)$ and $(3.28)$ into $\left( 3.21\right)
$ for sufficiently large $r\in E_{7}\backslash E_{6}$, we obtain%
\begin{equation*}
r^{\mu (A_{l0})-\varepsilon }\leq \left( (n-1)(m+1)+2m\right) r^{\rho
+\varepsilon }+r^{\mu (A_{kp})+\varepsilon }+O(r^{\rho (f)-1+\varepsilon })
\end{equation*}%
\begin{equation}
+r^{\lambda \left( \frac{1}{A_{l0}}\right) +\varepsilon }+O(\log r).
\tag{3.30}
\end{equation}%
We may choose sufficiently small $\varepsilon $ satisfying
\begin{equation*}
0<3\varepsilon <\min \left\{ \mu (A_{l0})-\rho ,\mu (A_{l0})-\lambda \left(
\frac{1}{A_{l0}}\right) \right\} ,
\end{equation*}%
from $\left( 3.30\right) $ that for $r\in E_{7}\backslash E_{6}$, $%
r\rightarrow +\infty $
\begin{equation*}
r^{\mu (A_{l0})-2\varepsilon }\leq r^{\rho (f)-1+\varepsilon },
\end{equation*}%
that means, $\mu (A_{l0})\leq \rho (f)-1+3\varepsilon $, since $\varepsilon
>0$ is arbitrary, then $\rho (f)\geq \mu (A_{l0})+1.$

\quad

\noindent \textbf{Case }$\left( \text{ii}\right) $ We suppose that $\max
\{\rho (A_{ij}),\rho (F),(i,j)\neq (l,0),(k,p)\}=\beta <\mu (A_{l0})=\mu
(A_{kp}),$ $\underline{\tau }(A_{l0})>\underline{\tau }(A_{kp}).$

\quad

\noindent For $S\in \left\{ F,A_{ij}:(i,j)\neq (l,0),(k,p)\right\} ,$ by the
definition of $\rho (S),$ for any given $\varepsilon (>0)$ and sufficiently
large $r$, we obtain
\begin{equation}
T(r,S)\leq r^{\rho (S)+\varepsilon }\leq r^{\beta +\varepsilon }.  \tag{3.31}
\end{equation}%
From the definition of $\underline{\tau }(A_{l0})$, for sufficiently small $%
\varepsilon >0$ and sufficiently large $r$, we have
\begin{equation}
T(r,A_{l0})\geq \left( \underline{\tau }(A_{l0})-\varepsilon \right) r^{\mu
(A_{l0})}.  \tag{3.32}
\end{equation}%
It follows from the definition of $\underline{\tau }(A_{kp})$ and Lemma
2.10, for any given $\varepsilon (>0)$, there exists a subset $E_{8}\subset
(1,+\infty )$ with infinite logarithmic measure such that for all $r\in
E_{8} $, we have
\begin{equation}
T(r,A_{kp})\leq \left( \underline{\tau }(A_{kp})+\varepsilon \right) r^{\mu
(A_{kp})}\leq \left( \underline{\tau }(A_{kp})+\varepsilon \right) r^{\mu
(A_{l0})}.  \tag{3.33}
\end{equation}%
By substituting $(3.22)-\left( 3.24\right) ,$ $\left( 3.26\right) $ and $%
\left( 3.31\right) -\left( 3.33\right) $ into $\left( 3.21\right) ,$ for
sufficiently large $r\in E_{8}\backslash E_{6}$, we obtain%
\begin{equation*}
\left( \underline{\tau }(A_{l0})-\varepsilon \right) r^{\mu (A_{l0})}\leq
\left( (n-1)(m+1)+2m\right) r^{\beta +\varepsilon }+\left( \underline{\tau }%
(A_{kp})+\varepsilon \right) r^{\mu (A_{l0})}
\end{equation*}%
\begin{equation}
+O(r^{\rho (f)-1+\varepsilon })+2\left( 2r\right) ^{\rho (f)+\varepsilon
}+r^{\beta +\varepsilon }+r^{\lambda \left( \frac{1}{A_{l0}}\right)
+\varepsilon }+O(\log r).  \tag{3.34}
\end{equation}%
Now, we may choose sufficiently small $\varepsilon $ satisfying%
\begin{equation*}
0<2\varepsilon <\min \left\{ \mu (A_{l0})-\beta ,\underline{\tau }(A_{l0})-%
\underline{\tau }(A_{kp}),\mu (A_{l0})-\lambda \left( \frac{1}{A_{l0}}%
\right) \right\} ,
\end{equation*}%
so from $(3.34)$ for sufficiently large $r\in E_{8}\backslash E_{6}$, we
have
\begin{equation*}
(\underline{\tau }(A_{l0})-\underline{\tau }(A_{kp})-2\varepsilon )r^{\mu
(A_{l0})-\varepsilon }\leq r^{\rho (f)+\varepsilon },
\end{equation*}%
that means, $\mu (A_{l0})\leq \rho (f)+2\varepsilon $, since $\varepsilon >0$
is arbitrary, then $\rho (f)\geq \mu (A_{l0}).$

\noindent Further, if \textit{\ }$F\equiv 0,$\textit{\ }then\textit{\ }by
substituting $(3.22)-(3.24)$ and $\left( 3.31\right) -\left( 3.33\right) $
into $\left( 3.21\right) ,$ for sufficiently large $r\in E_{8}\backslash
E_{6}$, we obtain%
\begin{equation*}
\left( \underline{\tau }(A_{l0})-\varepsilon \right) r^{\mu (A_{l0})}\leq
\left( (n-1)(m+1)+2m\right) r^{\beta +\varepsilon }+\left( \underline{\tau }%
(A_{kp})+\varepsilon \right) r^{\mu (A_{l0})}
\end{equation*}%
\begin{equation}
+O(r^{\rho (f)-1+\varepsilon })+r^{\lambda \left( \frac{1}{A_{l0}}\right)
+\varepsilon }+O(\log r).  \tag{3.35}
\end{equation}%
Now, we may choose sufficiently small $\varepsilon $ satisfying%
\begin{equation*}
0<2\varepsilon <\min \left\{ \mu (A_{l0})-\beta ,\underline{\tau }(A_{l0})-%
\underline{\tau }(A_{kp}),\mu (A_{l0})-\lambda \left( \frac{1}{A_{l0}}%
\right) \right\} .
\end{equation*}%
From $(3.35)$ for sufficiently large $r\in E_{8}\backslash E_{6}$, we get
\begin{equation*}
(\underline{\tau }(A_{l0})-\underline{\tau }(A_{kp})-2\varepsilon )r^{\mu
(A_{l0})-\varepsilon }\leq r^{\rho (f)-1+\varepsilon },
\end{equation*}%
that means, $\mu (A_{l0})\leq \rho (f)-1+2\varepsilon $, since $\varepsilon
>0$ is arbitrary, then $\rho (f)\geq \mu (A_{l0})+1.$

\noindent \textbf{Case }$\left( \text{iii}\right) $ We suppose that $\mu
(A_{l0})=\max \{\rho (F),\rho (A_{ij}):(i,j)\neq (l,0),(k,p)\}>\mu (A_{kp})$%
, and $\tau _{1}=\sum_{\rho (A_{ij})=\mu (A_{l0}),\,(i,j)\neq
(l,0),(k,p)}\tau (A_{ij})+\tau (F)<\underline{\tau }(A_{l0}).$ Then, there
exists a subset $J\subseteq \{0,1,\dots ,n\}\times \{0,1,\dots
,m\}\backslash \left\{ (l,0),(k,p)\right\} $ such that for all $(i,j)\in J,$
when $\rho (A_{ij})=\mu \left( A_{l0}\right) ,$ we have $\underset{(i,j)\in J%
}{\sum }\tau \left( A_{ij}\right) <\underline{\tau }\left( A_{l0}\right)
-\tau (F),$ and for $(i,j)\in \Pi =\{0,1,\dots ,n\}\times \{0,1,\dots
,m\}\backslash \left( J\cup \left\{ \,(l,0),(k,p)\right\} \right) $ we have $%
\rho \left( A_{ij}\right) <\mu \left( A_{l0}\right) .$ Hence, for any given $%
\varepsilon >0$ and sufficiently large $r,$ we get
\begin{equation}
T\left( r,A_{ij}\right) \leq \left\{
\begin{array}{c}
\left( \tau (A_{ij})+\varepsilon \right) r^{\mu (A_{l0})},\text{ if }%
(i,j)\in J, \\
r^{\rho (A_{ij})+\varepsilon }\leq r^{\mu (A_{l0})-\varepsilon },\text{ if }%
(i,j)\in \Pi%
\end{array}%
\right.  \tag{3.36}
\end{equation}%
and
\begin{equation}
T\left( r,F\right) \leq \left\{
\begin{array}{c}
\left( \tau (F)+\varepsilon \right) r^{\mu (A_{l0})},\text{ if }\rho (F)=\mu
(A_{l0}), \\
r^{\rho (F)+\varepsilon }\leq r^{\mu (A_{l0})-\varepsilon },\text{ if }\rho
(F)<\mu (A_{l0}).%
\end{array}%
\right.  \tag{3.37}
\end{equation}%
Then, by substituting $\left( 3.22\right) -\left( 3.24\right) ,$ $\left(
3.26\right) ,$ $\left( 3.28\right) ,$ $\left( 3.32\right) ,$ $\left(
3.36\right) $ and $\left( 3.37\right) $ into $\left( 3.21\right) $, for all $%
z$ satisfying $\left\vert z\right\vert =r\in E_{7}\backslash E_{6}$
sufficiently large $r$, we obtain%
\begin{equation*}
\left( \underline{\tau }(A_{l0})-\varepsilon \right) r^{\mu (A_{l0})}\leq
\underset{(i,j)\in J}{\sum }\left( \tau \left( A_{ij}\right) +\varepsilon
\right) r^{\mu \left( A_{l0}\right) }+\underset{(i,j)\in \Pi }{\sum }r^{\mu
(A_{l0})-\varepsilon }+r^{\mu (A_{kp})+\varepsilon }
\end{equation*}%
\begin{equation*}
+\left( \tau (F)+\varepsilon \right) r^{\mu (A_{l0})}+r^{\lambda \left(
\frac{1}{A_{l0}}\right) +\varepsilon }+2\left( 2r\right) ^{\rho
(f)+\varepsilon }+O\left( r^{\rho (f)-1+\varepsilon }\right) +O\left( \ln
r\right)
\end{equation*}%
\begin{equation*}
\leq \left( \tau _{1}+(nm+n+m)\varepsilon \right) r^{\mu \left(
A_{l0}\right) }+O\left( r^{\mu (A_{l0})-\varepsilon }\right) +r^{\mu
(A_{kp})+\varepsilon }
\end{equation*}%
\begin{equation}
+r^{\lambda \left( \frac{1}{A_{l0}}\right) +\varepsilon }+2\left( 2r\right)
^{\rho (f)+\varepsilon }+O\left( r^{\rho (f)-1+\varepsilon }\right) +O\left(
\log r\right) .  \tag{3.38}
\end{equation}%
Now, we choose sufficiently small $\varepsilon $ satisfying
\begin{equation*}
0<\varepsilon <\min \left\{ \frac{\mu (A_{l0})-\mu (A_{kp})}{2},\frac{%
\underline{\tau }(A_{l0})-\tau _{1}}{nm+n+m+1},\frac{\mu (A_{l0})-\lambda
\left( \frac{1}{A_{l0}}\right) }{2}\right\} ,
\end{equation*}%
then from $\left( 3.38\right) $ for sufficiently large $r\in E_{7}\backslash
E_{6}$, we get
\begin{equation*}
\left( \underline{\tau }(A_{l0})-\tau _{1}-\left( nm+n+m+1\right)
\varepsilon \right) r^{\mu (A_{l0})-\varepsilon }\leq r^{\rho
(f)+\varepsilon },
\end{equation*}%
that means, $\mu (A_{l0})\leq \rho (f)+2\varepsilon $, since $\varepsilon >0$
is arbitrary, then $\rho (f)\geq \mu (A_{l0}).$

\quad

\noindent Further, if \textit{\ }$F\equiv 0,$\textit{\ }then\textit{\ }by
substituting $\left( 3.22\right) -\left( 3.24\right) ,$ $\left( 3.28\right)
, $ $\left( 3.32\right) $ and $\left( 3.36\right) $ into $\left( 3.21\right)
$, for all $z$ satisfying $\left\vert z\right\vert =r\in E_{7}\backslash
E_{6} $ sufficiently large $r$, we obtain%
\begin{equation*}
\left( \underline{\tau }(A_{l0})-\varepsilon \right) r^{\mu (A_{l0})}\leq
\left( \tau _{1}+(nm+n+m-1)\varepsilon \right) r^{\mu \left( A_{l0}\right)
}+O\left( r^{\mu (A_{l0})-\varepsilon }\right)
\end{equation*}%
\begin{equation}
+r^{\mu (A_{kp})+\varepsilon }+r^{\lambda \left( \frac{1}{A_{l0}}\right)
+\varepsilon }+O\left( r^{\rho (f)-1+\varepsilon }\right) +O\left( \log
r\right) .  \tag{3.39}
\end{equation}%
Now, we choose sufficiently small $\varepsilon $ satisfying
\begin{equation*}
0<\varepsilon <\min \left\{ \frac{\mu (A_{l0})-\mu (A_{kp})}{2},\frac{%
\underline{\tau }(A_{l0})-\tau _{1}}{nm+n+m},\frac{\mu (A_{l0})-\lambda
\left( \frac{1}{A_{l0}}\right) }{2}\right\} ,
\end{equation*}%
then from $\left( 3.39\right) $ for sufficiently large $r\in E_{7}\backslash
E_{6}$, we get
\begin{equation*}
\left( \underline{\tau }(A_{l0})-\tau _{1}-\left( nm+n+m\right) \varepsilon
\right) r^{\mu (A_{l0})-\varepsilon }\leq r^{\rho (f)-1+\varepsilon },
\end{equation*}%
that means, $\mu (A_{l0})\leq \rho (f)-1+2\varepsilon $, since $\varepsilon
>0$ is arbitrary, then $\rho (f)\geq \mu (A_{l0})+1.$

\quad

\noindent \textbf{Case }$\left( \text{iv}\right) $ We suppose that $\max
\{\rho (A_{j}),\rho (F),(i,j)\neq (l,0),(k,p)\}=\mu (A_{l0})=\mu (A_{kp})$
with $\tau _{1}+\underline{\tau }(A_{kp})<\underline{\tau }(A_{l0}).$

\noindent It follows by substituting $\left( 3.22\right) -\left( 3.24\right)
,$ $\left( 3.26\right) ,$ $\left( 3.32\right) ,$ $\left( 3.33\right) ,$ $%
\left( 3.36\right) $ and $\left( 3.37\right) $ into $\left( 3.21\right) $,
for all sufficiently large $r\in E_{8}\backslash E_{6},$ we have%
\begin{equation*}
\left( \underline{\tau }(A_{l0})-\varepsilon \right) r^{\mu (A_{l0})}\leq
\left( \tau _{1}+(nm+n+m)\varepsilon \right) r^{\mu \left( A_{l0}\right)
}+O\left( r^{\mu (A_{l0})-\varepsilon }\right)
\end{equation*}%
\begin{equation}
+\left( \underline{\tau }(A_{kp})+\varepsilon \right) r^{\mu
(A_{l0})}+r^{\lambda \left( \frac{1}{A_{l0}}\right) +\varepsilon }+2\left(
2r\right) ^{\rho (f)+\varepsilon }+O\left( r^{\rho (f)-1+\varepsilon
}\right) +O\left( \log r\right) .  \tag{3.40}
\end{equation}%
Now, we choose sufficiently small $\varepsilon $ satisfying
\begin{equation*}
0<\varepsilon <\min \left\{ \frac{\underline{\tau }(A_{l0})-\tau _{1}-%
\underline{\tau }(A_{kp})}{nm+n+m+2},\frac{\mu (A_{l0})-\lambda \left( \frac{%
1}{A_{l0}}\right) }{2}\right\} ,
\end{equation*}%
then from $\left( 3.40\right) $ for sufficiently large $r\in E_{8}\backslash
E_{6}$, we get
\begin{equation*}
\left( \underline{\tau }(A_{l0})-\tau _{1}-\underline{\tau }(A_{kp})-\left(
nm+n+m+2\right) \varepsilon \right) r^{\mu (A_{l0})-\varepsilon }\leq
r^{\rho (f)+\varepsilon },
\end{equation*}%
that means, $\mu (A_{l0})\leq \rho (f)+2\varepsilon $, since $\varepsilon >0$
is arbitrary, then $\rho (f)\geq \mu (A_{l0}).$

\quad

\noindent Further, if \textit{\ }$F\equiv 0,$\textit{\ }then by It follows
by substituting $\left( 3.22\right) -\left( 3.24\right) ,$ $\left(
3.32\right) ,$ $\left( 3.33\right) $ and $\left( 3.36\right) $, for all
sufficiently large $r\in E_{8}\backslash E_{6},$ we have%
\begin{equation*}
\left( \underline{\tau }(A_{l0})-\varepsilon \right) r^{\mu (A_{l0})}\leq
\left( \tau _{1}+(nm+n+m-1)\varepsilon \right) r^{\mu \left( A_{l0}\right)
}+O\left( r^{\mu (A_{l0})-\varepsilon }\right)
\end{equation*}%
\begin{equation}
+\left( \underline{\tau }(A_{kp})+\varepsilon \right) r^{\mu
(A_{l0})}+r^{\lambda \left( \frac{1}{A_{l0}}\right) +\varepsilon }+O\left(
r^{\rho (f)-1+\varepsilon }\right) +O\left( \log r\right) .  \tag{3.41}
\end{equation}%
Now, we choose sufficiently small $\varepsilon $ satisfying
\begin{equation*}
0<\varepsilon <\min \left\{ \frac{\underline{\tau }(A_{l0})-\tau _{1}-%
\underline{\tau }(A_{kp})}{nm+n+m+1},\frac{\mu (A_{l0})-\lambda \left( \frac{%
1}{A_{l0}}\right) }{2}\right\} ,
\end{equation*}%
then from $\left( 3.41\right) $ for sufficiently large $r\in E_{8}\backslash
E_{6}$, we get
\begin{equation*}
\left( \underline{\tau }(A_{l0})-\tau _{1}-\underline{\tau }(A_{kp})-\left(
nm+n+m+1\right) \varepsilon \right) r^{\mu (A_{l0})-\varepsilon }\leq
r^{\rho (f)-1+\varepsilon },
\end{equation*}%
so that means, $\mu (A_{l0})\leq \rho (f)-1+2\varepsilon $, since $%
\varepsilon >0$ is arbitrary, then $\rho (f)\geq \mu (A_{l0})+1.$ The proof
of Theorem 1.2 is complete.

\section{\textbf{Examples}}

\noindent \textbf{Example 4.1} Consider the non-homogeneous
differential-difference equation with entire coefficients
\begin{equation}
A_{02}(z)f^{\prime \prime }(z)+A_{11}(z)f^{\prime }(z+1)+A_{01}(z)f^{\prime
}(z)+A_{10}(z)f(z+1)+A_{00}(z)f(z)=F\left( z\right) .  \tag{4.1}
\end{equation}%
\textbf{Case 1. }$\max \{\rho (F),\mu (A_{kp}),\rho (A_{ij}):(i,j)\neq
(l,0),(k,p)\}<\mu (A_{l0})$. In $\left( 4.1\right) $ for
\begin{equation*}
A_{00}(z)=\pi ^{2}+2\pi ^{4}z^{2},\quad A_{10}(z)=e^{-\pi ^{2}z^{2}-\pi
^{2}},\quad A_{01}(z)=2\pi ^{2}(z+1)e^{2\pi ^{2}z+\pi ^{2}},
\end{equation*}%
\begin{equation*}
A_{11}(z)=-2\pi ^{2}z,\quad A_{02}(z)=-\frac{1}{2},\quad F(z)=e^{2\pi ^{2}z},
\end{equation*}%
we have
\begin{equation*}
\max \{\rho (F),\mu (A_{11}),\rho (A_{ij}):(i,j)\neq (1,0),(1,1)\}=1<\mu
(A_{10})=2.
\end{equation*}%
We see that the conditions of Theorem 1.1 are verified. The function
\begin{equation*}
f\left( z\right) =e^{\pi ^{2}z^{2}}
\end{equation*}%
is a solution of equation $\left( 4.1\right) $ which satisfies $\rho
(f)=2\geq \mu (A_{10})=2$.

\quad

\noindent \textbf{Case 2.} $\max \{\rho (F),\rho (A_{ij}):(i,j)\neq
(l,0),(k,p)\}<\mu (A_{l0})=\mu (A_{kp})$ with $\underline{\tau }_{M}(A_{l0})>%
\underline{\tau }_{M}(A_{kp})$. In $\left( 4.1\right) $ for
\begin{equation*}
A_{00}(z)=2\pi ^{2},\quad A_{10}(z)=2\pi ^{2}(z+1)e^{z^{2}}+e^{-\pi
^{2}z^{2}-\pi ^{2}},\quad A_{01}(z)=2\pi ^{2}z,
\end{equation*}%
\begin{equation*}
A_{11}(z)=-e^{z^{2}},\quad A_{02}(z)=-1,\quad F(z)=e^{2\pi ^{2}z},
\end{equation*}%
we get $\max \{\rho (F),\rho (A_{ij}):(i,j)\neq (1,0),(1,1)\}=1<\mu
(A_{10})=\mu (A_{11})=2$ and $\underline{\tau }_{M}(A_{10})=\pi ^{2}>%
\underline{\tau }_{M}(A_{11})=1.$ Hence, the conditions of Theorem 1.1 are
satisfied. The function
\begin{equation*}
f\left( z\right) =e^{\pi ^{2}z^{2}}
\end{equation*}%
is a solution of equation $\left( 4.1\right) $ and $f$ satisfies $\rho
(f)=2\geq \mu (A_{10})=2$.

\quad

\noindent \textbf{Case 3.} $\mu (A_{l0})=\max \{\rho (F),\rho
(A_{ij}):(i,j)\neq (l,0),(k,p)\}>\mu (A_{kp})$ with $\underline{\tau }%
_{M}(A_{l0})>\tau _{1}=\max \{\tau _{M}(S):\rho (S)=\mu (A_{l0}),S\in
\{F,A_{ij}:(i,j)\neq (l,0),(k,p)\}\}$. In $\left( 4.1\right) $ for
\begin{equation*}
A_{00}(z)=\pi ^{2}+2\pi ^{4}z^{2},\quad A_{10}(z)=e^{-\frac{4}{5}\pi
^{2}z^{2}-\pi ^{2}},\quad A_{01}(z)=2\pi ^{2}(z+1)e^{2\pi ^{2}z+\pi ^{2}},
\end{equation*}%
\begin{equation*}
A_{11}(z)=-2\pi ^{2}z,\quad A_{02}(z)=-\frac{1}{2},\quad F(z)=e^{\frac{1}{5}%
\pi ^{2}z^{2}+2\pi ^{2}z},
\end{equation*}%
we have $\mu (A_{10})=\max \{\rho (F),\rho (A_{ij}):(i,j)\neq
(1,0),(0,1)\}=2>\mu (A_{01})=1$ and $\underline{\tau }(A_{10})=\frac{4\pi
^{2}}{5}>\tau _{1}=\tau _{M}(F)=\frac{\pi ^{2}}{5}.$ Obviously, the
conditions of Theorem 1.1 are verified. The function
\begin{equation*}
f\left( z\right) =e^{\pi ^{2}z^{2}}
\end{equation*}%
is a solution of equation $\left( 4.1\right) $ which satisfies $\rho
(f)=2\geq \mu (A_{10})=2$.

\quad

\noindent \textbf{Case 4.} $\mu (A_{l0})=\mu (A_{kp})=\max \{\rho (F),\rho
(A_{ij}):(i,j)\neq (l,0),(k,0)\},$ and $\underline{\tau }_{M}(A_{10})>\max
\{\tau _{1},\underline{\tau }_{M}(A_{kp})\}.$ In $\left( 4.1\right) $ for
\begin{equation*}
A_{00}(z)=2\pi ^{2},\quad A_{10}(z)=2\pi ^{2}(z+1)e^{z^{2}}+e^{-\frac{4}{5}%
\pi ^{2}z^{2}-\pi ^{2}},\quad A_{01}(z)=2\pi ^{2}z,
\end{equation*}%
\begin{equation*}
A_{11}(z)=-e^{z^{2}},\quad A_{02}(z)=-1,\quad F(z)=e^{\frac{1}{5}\pi
^{2}z^{2}+2\pi ^{2}z},
\end{equation*}%
we get $\mu (A_{10})=\mu (A_{11})=\max \{\rho (F),\rho (A_{ij}):(i,j)\neq
(1,0),(1,1)\}=2$ and $\underline{\tau }_{M}(A_{10})=\frac{4\pi ^{2}}{5}>\max
\{\tau _{1},\underline{\tau }_{M}(A_{11})\}=\max \{\tau _{M}(F),\underline{%
\tau }_{M}(A_{11})\}=\frac{\pi ^{2}}{5}.$ We see that the conditions of
Theorem 1.1 are satisfied. The function
\begin{equation*}
f\left( z\right) =e^{\pi ^{2}z^{2}}
\end{equation*}%
is a solution of equation $\left( 4.1\right) $ which satisfies $\rho
(f)=2\geq \mu (A_{10})=2$.

\quad

\noindent \textbf{Example 4.2} Consider the homogeneous
differential-difference equation with entire coefficients%
\begin{equation}
A_{11}(z)g^{\prime }(z-1)+A_{20}(z)g(z+3)+A_{00}(z)g(z)=0.  \tag{4.2}
\end{equation}%
\textbf{Case 1. }$\max \{\mu (A_{kp}),\rho (A_{ij}):(i,j)\neq
(l,0),(k,p)\}<\mu (A_{l0})$. In $\left( 4.2\right) ,$ for
\begin{equation*}
A_{00}(z)=1,\quad A_{20}(z)=\left( 4\pi i(1-z)-e^{4\pi iz}\right) e^{-16\pi
iz},\quad A_{11}(z)=1,
\end{equation*}%
we have
\begin{equation*}
\max \{\mu (A_{11}),\rho (A_{ij}):(i,j)\neq (2,0),(1,1)\}=0<\mu (A_{20})=1.
\end{equation*}%
So, the conditions of Theorem 1.1 are satisfied. The meromorphic function
\begin{equation*}
g\left( z\right) =e^{2\pi iz^{2}}
\end{equation*}%
is a solution of equation $\left( 4.2\right) $ and $g$ satisfies $\rho
(g)=2\geq \mu (A_{20})+1=2.$

\quad

\noindent \textbf{Case 2.} $\max \{\rho (A_{ij}):(i,j)\neq (l,0),(k,p)\}<\mu
(A_{l0})=\mu (A_{kp}),$ with $\underline{\tau }_{M}(A_{l0})>\underline{\tau }%
_{M}(A_{kp})$. In $\left( 4.2\right) ,$ for
\begin{equation*}
A_{00}(z)=1,\quad A_{20}(z)=\left( 4\pi i(1-z)-e^{2\pi iz}\right) e^{-14\pi
iz},\quad A_{11}(z)=e^{2\pi iz},
\end{equation*}%
we get $\mu (A_{20})=\mu (A_{11})=1>\max \{\rho (A_{ij}):(i,j)\neq
(2,0),(1,1)\}=0$ and $\underline{\tau }_{M}(A_{20})=14\pi >\underline{\tau }%
_{M}(A_{11})=2\pi .$ Obviously, the conditions of Theorem 1.1 are verified.
The meromorphic function
\begin{equation*}
g\left( z\right) =e^{2\pi iz^{2}}
\end{equation*}%
is a solution of equation $\left( 4.2\right) $ and $g$ satisfies $\rho
(g)=2\geq \mu (A_{20})+1=2.$

\quad

\noindent \textbf{Case 3.} $\mu (A_{l0})=\mu (A_{kp})=\max \{\rho
(A_{ij}):(i,j)\neq (l,0),(k,p)\}$ with $\underline{\tau }_{M}(A_{l0})>\tau
_{1}=\max \{\tau _{M}(A_{ij}):\rho (A_{ij})=\mu (A_{l0}),(i,j)\neq
(l,0),(k,p)\}$. In $\left( 4.2\right) ,$ for
\begin{equation*}
A_{00}(z)=e^{-2\pi iz},\quad A_{20}(z)=\left( 4\pi i(1-z)-1\right) e^{-14\pi
iz},\quad A_{11}(z)=e^{2\pi iz},
\end{equation*}%
we have $\mu (A_{20})=\mu (A_{11})=\max \{\rho (A_{ij}):(i,j)\neq
(2,0),(1,1)\}=1$ and $\underline{\tau }_{M}(A_{20})=14\pi >\max \{\tau
_{M}(A_{00}),\underline{\tau }_{M}(A_{11})\}=2\pi .$ It is clear that the
conditions of Theorem 1.1 are satisfied. The meromorphic function
\begin{equation*}
g\left( z\right) =e^{2\pi iz^{2}}
\end{equation*}%
is a solution of equation $\left( 4.2\right) $ and $g$ satisfies $\rho
(g)=2\geq \mu (A_{20})+1=2.$

\quad

\noindent \textbf{Example 4.3} Consider the non-homogeneous
differential-difference equation with meromorphic coefficients
\begin{equation}
A_{11}(z)f^{\prime }(z-1)+A_{01}(z)f^{\prime
}(z)+A_{20}(z)f(z+1)+A_{10}(z)f(z-1)=F\left( z\right) .  \tag{4.3}
\end{equation}%
\textbf{Case 1. }$\max \{\rho (F),\mu (A_{kp}),\rho (A_{ij}):(i,j)\neq
(l,0),(k,p)\}<\mu (A_{l0})$. In $\left( 4.3\right) $ for%
\begin{equation*}
A_{10}(z)=e^{-\pi ^{3}z^{3}+3\pi ^{3}z^{2}-3\pi ^{3}z+\pi ^{3}},\quad
A_{20}(z)=3\pi ^{3}\left( 2z-1\right) e^{-3\pi ^{3}z^{2}-3\pi ^{3}z-\pi
^{3}},
\end{equation*}%
\begin{equation*}
A_{01}(z)=-1,\quad A_{11}(z)=e^{3\pi ^{3}z^{2}-3\pi ^{3}z+\pi ^{3}},\quad
F(z)=\tan (\pi z),
\end{equation*}%
we have
\begin{equation*}
\max \{\rho (F),\mu (A_{11}),\rho (A_{ij}):(i,j)\neq (1,0),(1,1)\}=2<\mu
(A_{10})=3,
\end{equation*}%
\begin{equation*}
\lambda \left( \frac{1}{A_{10}}\right) =0<\mu \left( A_{10}\right) =3.
\end{equation*}%
It is easy to see that the conditions of Theorem 1.2 are verified. The
function
\begin{equation*}
f\left( z\right) =e^{\pi ^{3}z^{3}}\tan (\pi z)
\end{equation*}%
is a solution of equation $\left( 4.3\right) $ which satisfies $\rho
(f)=3\geq \mu (A_{10})=3$.

\quad

\noindent \textbf{Case 2.} $\max \{\rho (F),\rho (A_{ij}):(i,j)\neq
(l,0),(k,p)\}<\mu (A_{l0})=\mu (A_{kp})$ with $\underline{\tau }(A_{l0})>%
\underline{\tau }(A_{kp})$. In $\left( 4.3\right) $ for%
\begin{equation*}
A_{10}(z)=e^{-\pi ^{3}z^{3}+3\pi ^{3}z^{2}-3\pi ^{3}z+\pi ^{3}}
\end{equation*}%
\begin{equation*}
+\left( 3\pi ^{3}\left( 2z-1\right) -3z^{2}\pi ^{3}-\pi \tan (\pi z)+\frac{%
\pi }{\tan (\pi z)}\right) e^{-z^{3}},
\end{equation*}%
\begin{equation*}
A_{20}(z)=3\pi ^{3}z^{2}+\pi \tan (\pi z)+\frac{\pi }{\tan (\pi z)},\quad
A_{01}(z)=-e^{3\pi ^{3}z^{2}+3\pi ^{3}z+\pi ^{3}},
\end{equation*}%
\begin{equation*}
A_{11}(z)=e^{-z^{3}},\quad F(z)=\tan (\pi z),
\end{equation*}%
we get
\begin{equation*}
\max \{\rho (F),\rho (A_{ij}):(i,j)\neq (1,0),(1,1)\}=2<\mu (A_{10})=\mu
(A_{11})=3,
\end{equation*}%
\begin{equation*}
\lambda \left( \frac{1}{A_{10}}\right) =1<\mu \left( A_{10}\right) =3
\end{equation*}%
and
\begin{equation*}
\underline{\tau }(A_{10})=\pi ^{2}>\underline{\tau }(A_{11})=\frac{1}{\pi }.
\end{equation*}%
Hence, the conditions of Theorem 1.2 are satisfied. The function
\begin{equation*}
f\left( z\right) =e^{\pi ^{3}z^{3}}
\end{equation*}%
is a solution of equation $\left( 4.3\right) $ and $f$ satisfies $\rho
(f)=3\geq \mu (A_{10})=3$.

\quad

\noindent \textbf{Case 3.} $\mu (A_{l0})=\max \{\rho (F),\rho
(A_{ij}):(i,j)\neq (l,0),(k,p)\}>\mu (A_{kp})$ with $\underline{\tau }%
(A_{l0})>\tau _{1}=\sum_{\rho (A_{ij})=\mu (A_{l0}),\,(i,j)\neq
(l,0),(k,p)}\tau (A_{ij})+\tau (F)$. In $\left( 4.3\right) $ for
\begin{equation*}
A_{10}(z)=e^{-2\pi ^{3}z^{3}+3\pi ^{3}z^{2}-3\pi ^{3}z+\pi ^{3}},\quad
A_{20}(z)=3\pi ^{3}\left( 2z-1\right) e^{-3\pi ^{3}z^{2}-3\pi ^{3}z-\pi
^{3}},
\end{equation*}%
\begin{equation*}
A_{01}(z)=-1,\quad A_{11}(z)=e^{3\pi ^{3}z^{2}-3\pi ^{3}z+\pi ^{3}},\quad
F(z)=\frac{\tan (\pi z)}{e^{\pi ^{3}z^{3}}},
\end{equation*}%
we have
\begin{equation*}
\mu (A_{10})=\max \{\rho (F),\rho (A_{ij}):(i,j)\neq (1,0),(1,1)\}=3>\mu
(A_{11})=2,
\end{equation*}%
\begin{equation*}
\lambda \left( \frac{1}{A_{10}}\right) =0<\mu \left( A_{10}\right) =3
\end{equation*}%
and $\underline{\tau }(A_{10})=2\pi ^{2}>\tau _{1}=\tau (F)=\pi ^{2}.$ We
can see that the conditions of Theorem 1.2 are verified. The function
\begin{equation*}
f\left( z\right) =e^{\pi ^{3}z^{3}}\tan (\pi z)
\end{equation*}%
is a solution of equation $\left( 4.3\right) $ which satisfies $\rho
(f)=3\geq \mu (A_{10})=3$.

\quad

\noindent \textbf{Case 4.} $\mu (A_{l0})=\mu (A_{k0})=\max \{\rho (F),\rho
(A_{ij}),:(i,j)\neq (l,0),(k,0)\},$ and $\underline{\tau }(A_{10})>\tau _{1}+%
\underline{\tau }(A_{kp}).$ In $\left( 4.3\right) $ for%
\begin{equation*}
A_{10}(z)=e^{-2\pi ^{3}z^{3}+3\pi ^{3}z^{2}-3\pi ^{3}z+\pi ^{3}},\quad
A_{20}(z)=3\pi ^{3}\left( 2z-1\right) e^{(\frac{\pi }{4}z)^{3}-3\pi
^{3}z^{2}-3\pi ^{3}z-\pi ^{3}},
\end{equation*}%
\begin{equation*}
A_{01}(z)=-e^{(\frac{\pi }{4}z)^{3}},\quad \ A_{11}(z)=e^{(\frac{\pi }{4}%
z)^{3}+3\pi ^{3}z^{2}-3\pi ^{3}z+\pi ^{3}},\quad F(z)=\frac{\tan (\pi z)}{%
e^{\pi ^{3}z^{3}}},
\end{equation*}%
we get
\begin{equation*}
\mu (A_{10})=\mu (A_{11})=\max \{\rho (F),\rho (A_{ij}):(i,j)\neq
(1,0),(1,1)\}=3,
\end{equation*}%
\begin{equation*}
\lambda \left( \frac{1}{A_{10}}\right) =0<\mu \left( A_{10}\right) =3
\end{equation*}%
and $\tau _{1}+\underline{\tau }(A_{11})=\tau (A_{01})+\tau (A_{20})+\tau
(F)+\underline{\tau }(A_{11})=\left( \frac{2}{4^{3}}+1\right) \pi ^{2}+\frac{%
\pi ^{2}}{4^{3}}=\frac{67}{64}\pi ^{2}<\underline{\tau }(A_{10})=2\pi ^{2}.$
Obviously, the conditions of Theorem 1.2 are satisfied. The function
\begin{equation*}
f\left( z\right) =e^{\pi ^{3}z^{3}}\tan (\pi z)
\end{equation*}%
is a solution of equation $\left( 4.3\right) $ which satisfies $\rho
(f)=3\geq \mu (A_{10})=3$.

\quad

\noindent \textbf{Example 4.4} Consider the homogeneous
(differential)-difference equation with meromorphic coefficients%
\begin{equation}
A_{11}(z)h^{\prime }(z+i\pi )+A_{20}(z)h(z+2i\pi )+A_{00}(z)h(z)=0.
\tag{4.4}
\end{equation}%
\textbf{Case 1. }$\max \{\mu (A_{kp}),\rho (A_{ij}):(i,j)\neq
(l,0),(k,p)\}<\mu (A_{l0})$. In $\left( 4.4\right) ,$ for
\begin{equation*}
A_{00}(z)=-1,\quad A_{20}(z)=e^{12\pi z^{2}+24\pi ^{2}iz-16\pi ^{3}}-e^{6\pi
z^{2}+18\pi ^{2}iz-14\pi ^{3}},
\end{equation*}%
\begin{equation*}
A_{11}(z)=\frac{\cos (2iz)}{6i(z+i\pi )^{2}\cos (2iz)+2i\sin (2iz)},
\end{equation*}%
we have $\max \{\mu (A_{11}),\rho (A_{ij}):(i,j)\neq (2,0),(1,1)\}=1<\mu
(A_{20})=2$ and
\begin{equation*}
\lambda \left( \frac{1}{A_{20}}\right) =0<\mu \left( A_{20}\right) =2.
\end{equation*}%
Obviously, the conditions of Theorem 1.2 are verified. The meromorphic
function
\begin{equation*}
h\left( z\right) =\frac{e^{2iz^{3}}}{\cos (2iz)}
\end{equation*}%
is a solution of equation $\left( 4.4\right) $ and $h$ satisfies $\rho
(h)=3\geq \mu (A_{20})+1=3.$

\noindent \textbf{Case 2.} $\max \{\rho (A_{ij}):(i,j)\neq (l,0),(k,p)\}<\mu
(A_{l0})=\mu (A_{kp})$ with $\underline{\tau }(A_{l0})>\underline{\tau }%
(A_{kp})$. In $\left( 4.4\right) ,$ for
\begin{equation*}
A_{00}(z)=1,\quad A_{20}(z)=-2e^{12\pi z^{2}+24\pi ^{2}iz-16\pi ^{3}},
\end{equation*}%
\begin{equation*}
A_{11}(z)=\frac{e^{6\pi z^{2}+6\pi ^{2}iz-2\pi ^{3}}\cos (2iz)}{6i(z+i\pi
)^{2}\cos (2iz)+2i\sin (2iz)},
\end{equation*}%
we get
\begin{equation*}
\mu (A_{20})=\mu (A_{11})=2>\max \{\rho (A_{ij}):(i,j)\neq
(2,0),(1,1)\}=\rho (A_{00})=0,
\end{equation*}%
\begin{equation*}
\lambda \left( \frac{1}{A_{20}}\right) =0<\mu \left( A_{20}\right) =2
\end{equation*}%
and
\begin{equation*}
\underline{\tau }(A_{20})=12>\underline{\tau }(A_{11})=6.
\end{equation*}%
It is clear that the conditions of Theorem 1.2 are satisfied. The
meromorphic function
\begin{equation*}
h\left( z\right) =\frac{e^{2iz^{3}}}{\cos (2iz)}
\end{equation*}%
is a solution of equation $\left( 4.4\right) $ and $h$ satisfies $\rho
(h)=3\geq \mu (A_{20})+1=3.$

\noindent \textbf{Case 3.} $\mu (A_{l0})=\mu (A_{kp})=\max \{\rho
(A_{ij}):(i,j)\neq (l,0),(k,p)\},$ with $\underline{\tau }%
(A_{l0})>\sum_{\rho (A_{ij})=\mu (A_{l0}),\,(i,j)\neq (l,0),(k,p)}\tau
(A_{ij})+\underline{\tau }(A_{kp})$. In $\left( 4.4\right) ,$ for
\begin{equation*}
A_{00}(z)=e^{\pi z^{2}},\quad A_{20}(z)=-2e^{13\pi z^{2}+24\pi ^{2}iz-16\pi
^{3}},
\end{equation*}%
\begin{equation*}
A_{11}(z)=\frac{e^{7\pi z^{2}+6\pi ^{2}iz-2\pi ^{3}}\cos (2iz)}{6i(z+i\pi
)^{2}\cos (2iz)+2i\sin (2iz)},
\end{equation*}%
we have
\begin{equation*}
\mu (A_{20})=\mu (A_{11})=\max \{\rho (A_{ij}):(i,j)\neq (2,0),(1,1)\}=\rho
(A_{00})=2,
\end{equation*}%
\begin{equation*}
\lambda \left( \frac{1}{A_{20}}\right) =0<\mu \left( A_{20}\right) =2
\end{equation*}%
and
\begin{equation*}
\underline{\tau }(A_{20})=13>\tau (A_{00})+\underline{\tau }(A_{11})=1+7=8.
\end{equation*}%
It is easy to see that the conditions of Theorem 1.2 are verified. The
meromorphic function
\begin{equation*}
h\left( z\right) =\frac{e^{2iz^{3}}}{\cos (2iz)}
\end{equation*}%
is a solution of equation $\left( 4.4\right) $ and $h$ satisfies $\rho
(h)=3\geq \mu (A_{20})+1=3.$

\quad

\noindent \textbf{Acknowledgements.} This paper was supported by the
Directorate-General for Scientific Research and Technological Development
(DGRSDT).

\begin{center}
{\Large References}
\end{center}

\noindent $\left[ 1\right] $ B. Bela\"{\i}di and Y. Benkarouba, \textit{Some
properties of meromorphic solutions of higher order linear difference
equation}s. Ser. A: Appl. Math. Inform. and Mech. vol. 11, 2 (2019), 75--95.

\noindent $\left[ 2\right] $ B. Bela\"{\i}di and R. Bellaama, \textit{%
Meromorphic solutions of higher order non-homogeneous linear difference
equations}. Bulletin of the Transilvania University of Brasov, Series III:
Mathematics, Informatics, Physics. Vol. 13, Issue 2 (2020), \ 433--450.

\noindent $\left[ 3\right] $ B. Bela\"{\i}di and R. Bellaama, \textit{Study
of the growth properties of meromorphic solutions of higher-order linear
difference equations}. Arab. J. Math. (2021).
https://doi.org/10.1007/s40065-021-00324-2

\noindent $\lbrack 4]$ Z. X. Chen, \textit{The zero, pole and orders of
meromorphic solutions of differential equations with meromorphic coefficients%
}. Kodai Math. J. 19 (1996), no. 3, 341--354.

\noindent $\lbrack 5]$ Z. X. Chen, \textit{Complex differences and
difference equations}, Mathematics Monograph Series 29. Science Press,
Beijing (2014).

\noindent $\left[ 6\right] $ Z. Chen and X. M. Zheng, \textit{Growth of
meromorphic solutions of general complex linear differential-difference
equation}. Acta Univ. Apulensis Math. Inform. No. 56 (2018), 1--12.

\noindent $\left[ 7\right] $ Y. M. Chiang and S. J. Feng, \textit{On the
Nevanlinna characteristic of }$f\left( z+\eta \right) $ \textit{and
difference equations in the complex plane. }Ramanujan J. 16 (2008), no. 1,
105--129.

\noindent $\left[ 8\right] $ A. Goldberg and I. Ostrovskii, \textit{Value
distribution of meromorphic functions}. Transl. Math. Monogr., vol. 236,
Amer. Math. Soc., Providence RI, 2008.

\noindent $\left[ 9\right] $ G. G. Gundersen, \textit{Estimates for the
logarithmic derivative of a meromorphic function, plus similar estimates}.\
J. London Math. Soc. (2) 37 (1988), no. 1, 88--104.

\noindent $\left[ 10\right] $ R. G. Halburd and R. J. Korhonen, \textit{%
Difference analogue of the lemma on the logarithmic derivative with
applications to difference equations. }J. Math. Anal. Appl. 314 (2006)%
\textit{, }no. 2, 477--487.

\noindent $\left[ 11\right] $ R. G. Halburd and R. J. Korhonen, \textit{%
Nevanlinna theory for the difference operator}. Ann. Acad. Sci. Fenn. Math.
31 (2006), no. 2, 463--478.

\noindent $\left[ 12\right] $ W. K. Hayman, \textit{Meromorphic functions}.
Oxford Mathematical Monographs Clarendon Press, Oxford 1964.

\noindent $\left[ 13\right] \ $H. Hu and X. M. Zheng, \textit{Growth of
solutions to linear differential equations with entire coefficients}.
Electron. J. Differential Equations 2012, No. 226, 15 pp.

\noindent $\left[ 14\right] $ I. Laine, \textit{Nevanlinna theory and
complex differential equations}. de Gruyter Studies in Mathematics, 15.
Walter de Gruyter \& Co., Berlin, 1993.

\noindent $\left[ 15\right] $ I. Laine and C. C. Yang, \textit{Clunie
theorems for difference and q-difference polynomials}. J. Lond. Math. Soc.
(2) 76 (2007), no. 3, 556-566.

\noindent $\left[ 16\right] $ Z. Latreuch and B. Bela\"{\i}di, \textit{%
Growth and oscillation of meromorphic solutions of linear difference
equations}. Mat. Vesnik 66 (2014), no. 2, 213--222.

\noindent $\left[ 17\right] $ H. F. Liu and Z. Q. Mao\textbf{, }\textit{On
the meromorphic solutions of some linear difference equations}. Adv.
Difference Equ. 2013, 2013:133, 1--12.

\noindent $\left[ 18\right] $ K. Liu and C. J. Song. \textit{Meromorphic
solutions of complex differential-difference equations}. Results Math. 72
(2017), no. 4, 1759--1771.

\noindent $\left[ 19\right] $ I. Q. Luo and X. M. Zheng\textbf{, }\textit{%
Growth of meromorphic solutions of some kind of complex linear difference
equation with entire or meromorphic coefficients}. Math. Appl. (Wuhan) 29
(2016), no. 4, 723--730.

\noindent $\left[ 20\right] $ A. G. Naftalevi\v{c}, \textit{Meromorphic
solutions of a differential-difference equation}. (Russian) Uspehi Mat. Nauk
16 1961 no. 3 (99), 191--196.

\noindent $\left[ 21\right] $ X. G. Qi and L. Z. Yang. \textit{A note on
meromorphic solutions of complex differential-difference equations}. Bull.
Korean Math. Soc. 56 (2019), no. 3, 597--607.

\noindent $\left[ 22\right] \ $S. Z. Wu and X. M. Zheng, \textit{Growth of
meromorphic solutions of complex linear differential-difference equations
with coefficients having the same order}. J. Math. Res. Appl. 34 (2014), no.
6, 683--695.\textit{\ }

\noindent $\left[ 23\right] \ $S. Z. Wu and X. M. Zheng, \textit{Growth of
solutions to some higher-order linear differential equations in }$%
\mathbb{C}
$\textit{\ and in unit disc} $\Delta $. (Chinese) Math. Appl. (Wuhan) 29
(2016), no. 1, 20--30.

\noindent $\left[ 24\right] $ C. C. Yang and H. X. Yi, \textit{Uniqueness
theory of meromorphic functions}. Mathematics and its Applications, 557.
Kluwer Academic Publishers Group, Dordrecht, 2003.

\noindent $\left[ 25\right] $ X. M. Zheng and J. Tu, \textit{Growth of
meromorphic solutions of linear difference equations}. J. Math. Anal. Appl.
384 (2011), no. 2, 349--356.

\noindent $\left[ 26\right] $ Y. P. Zhou and X. M. Zheng,\textbf{\ }\textit{%
Growth of meromorphic solutions to homogeneous and non-homogeneous linear
(differential-)difference equations with meromorphic coefficients}.
Electron. J. Differential Equations 2017, Paper No. 34, 15 pp.

\end{document}